\documentclass{amsart}
\RequirePackage{amssymb}
\usepackage{hyperref}
\calclayout

\newtheorem{theorem}{Theorem}[section]
\newtheorem{lemma}[theorem]{Lemma}

\newtheorem{corollary}[theorem]{Corollary}
\numberwithin{equation}{section}
\theoremstyle{definition}
\newtheorem{remark}[theorem]{Remark}
\newcommand{\C}{\mathbb{C}}
\newcommand{\N}{\mathbb{N}}
\newcommand{\T}{\mathbb{T}}
\newcommand{\Z}{\mathbb{Z}}
\newcommand{\R}{\mathbb{R}}
\newcommand{\cA}{\mathcal{A}}
\newcommand{\cB}{\mathcal{B}}
\newcommand{\cC}{\mathcal{C}}
\newcommand{\cJ}{\mathcal{J}}
\newcommand{\cL}{\mathcal{L}}
\newcommand{\cM}{\mathcal{M}}
\newcommand{\cP}{\mathcal{P}}
\newcommand{\cQ}{\mathcal{Q}}
\newcommand{\cR}{\mathcal{R}}
\newcommand{\diag}{\mathrm{diag}\,}
\newcommand{\im}{\mathrm{Im}\,}
\newcommand{\ind}{\mathrm{Ind}\,}
\newcommand{\Ker}{\mathrm{Ker}\,}
\newcommand{\tr}{\mathrm{trace}}

\begin{document}
\title[Krein algebras and Toeplitz determinants]
{Generalized Krein algebras and\\ asymptotics of Toeplitz determinants}

\author{A. B\"ottcher}
\address{%
Fakult\"at f\"ur Mathematik,
Technische Universit\"at Chemnitz,
D-09107, Chemnitz, Germany}
\email{aboettch@mathematik.tu-chemnitz.de}

\author{A. Karlovich}
\address{%
Departamento de Matem\'atica,
Instituto Superior T\'ecnico,
Av. Rovisco Pais 1,
1049-001, Lisbon, Portugal}
\email{akarlov@math.ist.utl.pt}
\thanks{The second author is partially supported by F.C.T. (Portugal) grants
SFRH/BPD/11619/2002 and FCT/FEDER/POCTI/MAT/59972/2004}

\author{B. Silbermann}
\address{%
Fakult\"at f\"ur Mathematik,
Technische Universit\"at Chemnitz,
D-09107, Chemnitz, Germany}
\email{silbermn@mathematik.tu-chemnitz.de}

\subjclass[2000]{Primary 47B35; Secondary 15A15, 47B10}
\date{\today}
\keywords{Toeplitz matrix, generalized Krein algebra, Szeg\H{o}-Widom limit theorem,
Toeplitz operator, Hankel operator, Wiener-Hopf factorization, Besov space,
H\"older space, Schatten-von Neumann class, regularized operator determinant}
\dedicatory{This paper is dedicated to the centenary of Mark Krein (1907--1989).}
\begin{abstract}
We give a survey on generalized Krein algebras $K_{p,q}^{\alpha,\beta}$
and their applications to Toeplitz determinants.
Our methods originated in a paper by Mark Krein of 1966, where he showed that
$K_{2,2}^{1/2,1/2}$ is a Banach algebra. Subsequently, Widom
proved the strong Szeg\H{o} limit theorem for block Toeplitz
determinants with symbols in $(K_{2,2}^{1/2,1/2})_{N\times N}$
and later two of the authors studied symbols in the generalized Krein
algebras $(K_{p,q}^{\alpha,\beta})_{N\times N}$, where
$\lambda:=1/p+1/q=\alpha+\beta$ and $\lambda=1$. We here extend these results
to $0<\lambda<1$. The entire paper is based on fundamental work by Mark Krein,
ranging from operator ideals through Toeplitz operators up to Wiener-Hopf factorization.
\end{abstract}

\maketitle
\section{Introduction and main results}
\subsection{Krein algebras}
Suppose that $\cA$ is a complex Banach space with norm $\|\cdot\|_\cA$ which
is also an algebra with unit $e$ over the field of complex numbers and
$\|e\|_\cA\ne 0$. If the multiplication in $\cA$ is continuous,
then $\cA$ is called a \textit{unital Banach algebra}. Such an algebra
can be equipped with a new norm $\|\cdot\|$ which is equivalent
to $\|\cdot\|_\cA$ and satisfies
\begin{equation}\label{eq:Banach-norm}
\|e\|=1,\quad\quad
\|ab\|\le\|a\|\,\|b\|\quad\mbox{for all}\quad a,b\in\cA.
\end{equation}
Each norm satisfying (\ref{eq:Banach-norm}) is called a \textit{Banach algebra norm}.

Let $\T$ be the unit circle and, for $1\le p\le\infty$, let $L^p:=L^p(\T)$ and
$H^p:=H^p(\T)$ be the standard Lebesgue and Hardy spaces. Denote by $\{a_k\}_{k\in\Z}$
the sequence of the Fourier coefficients of a function $a\in L^1$,
\[
a_k=\frac{1}{2\pi}\int_0^{2\pi}a(e^{i\theta})e^{-ik\theta}d\theta
\quad (k\in\Z).
\]

It was Mark Krein \cite{Kr66} who first discovered that the set of all functions
$a\in L^\infty$ satisfying $\sum_{k\in\Z}|a_k|^2|k|<\infty$ forms a Banach algebra.
This algebra is called the \textit{Krein algebra} and is denoted by $K_{2,2}^{1/2,1/2}$.

Now we give an equivalent definition of the Krein algebra. For $k\in\Z$, let
$\chi_k(t):=t^k$, where $t\in\T$. Let $\sum_{k\in\Z}a_k\chi_k$ be the Fourier
series of a function $a\in L^1$. The Riesz projections $P$ and $Q$ are defined
formally by
\[
P:\sum_{k\in\Z}a_k\chi_k\mapsto \sum_{k\ge 0}a_k\chi_k,
\quad
Q:\sum_{k\in\Z}a_k\chi_k\mapsto \sum_{k<0}a_k\chi_k.
\]
These operators are bounded on $L^2$. Besides the projections $P$ and
$Q$, we now need the so-called flip operator $J$. This is the isometric
operator acting on $L^p$, $1<p<\infty$, by $(Jf)(t)=(1/t)\widetilde{f}(t)$, where
$\widetilde{f}(t):=f(1/t)$. For $a\in L^\infty$, we define the Hankel operators
$H(a)$ and $H(\widetilde{a})$ by
\begin{equation}\label{eq:def-Hankel}
H(a):H^2\to H^2,
\quad f\mapsto PM(a)QJf;
\quad\quad
H(\widetilde{a}):H^2\to H^2,
\quad
f\mapsto JQM(a)Pf,
\end{equation}
where $M(a)f=af$ is the multiplication operator by $a$. It is well known
and easy to see that if $a\in L^\infty$, then both $H(a)$ and $H(\widetilde{a})$
are Hilbert-Schmidt if and only if $a\in K_{2,2}^{1/2,1/2}$. Thus,
\[
K_{2,2}^{1/2,1/2}=\{a\in L^\infty: H(a), H(\widetilde{a})\ \mbox{are Hilbert-Schmidt}\}.
\]

The first natural generalization of the classical Krein algebra $K_{2,2}^{1/2,1/2}$
consists in replacing the ideal of Hilbert-Schmidt operators by the Schatten-von Neumann
ideals $\cC_p(H^2)$, $1\le p\le\infty$ (see Section~\ref{sec:Schatten}).

Let $1\le p\le\infty$ and $1\le q\le\infty$ (but not necessarily $1/p+1/q=1$).
We put
\begin{equation}\label{eq:def-Krein}
\begin{split}
&
K_{p,0}^{1/p,0}
:=
\big\{a\in L^\infty :  H(\widetilde{a})\in\cC_p(H^2)\big\},
\\
&
K_{0,q}^{0,1/q}
:=
\big\{a\in L^\infty :  H(a)\in\cC_q(H^2)\big\},
\\
& K_{p,q}^{1/p,1/q} := K_{p,0}^{1/p,0}\cap K_{0,q}^{0,1/q}.
\end{split}
\end{equation}
It is clear that the sets (\ref{eq:def-Krein})
are linear spaces. We define norms by
\begin{equation}\label{eq:norm-Krein}
\begin{split}
&
\|a\| :=
\|a\|_{L^\infty}+\|H(\widetilde{a})\|_{\cC_p(H^2)}
\quad (a\in K_{p,0}^{1/p,0}),
\\
&
\|a\| :=
\|a\|_{L^\infty}+\|H(a)\|_{\cC_q(H^2)}
\quad (a\in K_{0,q}^{0,1/q}),
\\
&
\|a\| :=
\|a\|_{L^\infty}
+\|H(\widetilde{a})\|_{\cC_p(H^2)}
+\|H(a)\|_{\cC_q(H^2)}
\quad (a\in K_{p,q}^{1/p,1/q}).
\end{split}
\end{equation}

Let $A,B\subset L^\infty$ be two subsets. We put
$A+B=\{\varphi+\psi:\varphi\in A,\psi\in B\}$. Let
$\overline{H^p}$ be the set of functions $f$ in $L^p$ such that
the complex conjugate $\overline{f}$ belongs to $H^p$ and let
$C:=C(\T)$ be the set of continuous functions. It is well known
that $C+H^\infty$, and hence, $C+\overline{H^\infty}$ and
$QC:=(C+H^\infty)\cap (C+\overline{H^\infty})$ are closed
subalgebras of $L^\infty$ (see, e.g. \cite[Section~6.31]{D98} or
\cite[Chap.~1, Theorem~5.1]{Pe03}).

Hartman's theorem (see, e.g. \cite[Theorem~2.54]{BS06} or
\cite[Chap.~1, Theorem~5.5]{Pe03}) shows that
\[
K_{\infty,0}^{1/\infty,0} =C+H^\infty,
\quad
K_{0,\infty}^{0,1/\infty}=C+\overline{H^\infty},
\quad
K_{\infty,\infty}^{1/\infty,1/\infty}=QC,
\]
and since $\cC_p(H^2)\subset\cC_\infty(H^2)$, we have
\begin{equation}\label{eq:embedding}
K_{p,0}^{1/p,0}\subset C+H^\infty,
\quad
K_{0,q}^{0,1/q}\subset C+\overline{H^\infty},
\quad
K_{p,q}^{1/p,1/q}\subset QC.
\end{equation}

For a unital Banach algebra $\cA$, we denote by $G\cA$ its group of
invertible elements.
\begin{theorem}\label{th:Krein-invertibility}
Let $1\le p,q\le\infty$.
\begin{enumerate}
\item[(a)]
The sets {\rm (\ref{eq:def-Krein})} are Banach
algebras under the Banach algebra norms {\rm (\ref{eq:norm-Krein})}.
\item[(b)]
If $a\in K_{p,0}^{1/p,0}$, then
$a\in GK_{p,0}^{1/p,0}
\Longleftrightarrow
a\in G(C+H^\infty)$.

\item[(c)]
If $a\in K_{0,q}^{0,1/q}$, then
$a\in GK_{0,q}^{0,1/q}
\Longleftrightarrow
a\in G(C+\overline{H^\infty})$.

\item[(d)]
If $a\in K_{p,q}^{1/p,1/q}$, then
\[
a\in GK_{p,q}^{1/p,1/q}
\Longleftrightarrow
a\in  G(C+H^\infty)
\Longleftrightarrow
a\in G(C+\overline{H^\infty})
\Longleftrightarrow
a\in GL^\infty.
\]
\end{enumerate}
\end{theorem}
This theorem was established by Krein \cite{Kr66} for the case $p=q=2$.
We therefore call (\ref{eq:def-Krein}) \emph{Krein algebras}.
Theorem~\ref{th:Krein-invertibility} was proved in this form for the first time
in \cite[Sections~4.10--4.11]{BS83}, a complete proof is also given in
\cite[Theorem~10.9]{BS06}.
\subsection{Hankel operators in Schatten-von Neumann classes}
Of course, it is desirable to have equivalent definitions for $K_{p,q}^{1/p,1/q}$,
$K_{p,0}^{1/p,0}$, and $K_{0,q}^{0,1/q}$ in terms of functions and not operators.
As we already noted, for $K_{2,2}^{1/2,1/2}$ this is very easy. So, we need
effective criteria guaranteeing that $H(a)\in \cC_p(H^2)$ and $H(\widetilde{a})\in\cC_q(H^2)$
if $p\ne 2$ and $q\ne 2$.

Let $X$ be one of the spaces $L^p$, $1\le p<\infty$ or $C$. The moduli of
continuity of $f\in X$ are defined for $s\ge 0$ by
\[
\begin{split}
\omega_X^1(f,s)
&:=
\sup_{|h|\le s}\|f(e^{i(\cdot+h)})-f(e^{i\cdot})\|_X,
\\
\omega_X^2(f,s) &:= \sup_{|h|\le
s}\|f(e^{i(\cdot+h)})-2f(e^{i\cdot})+f(e^{i(\cdot-h)})\|_X.
\end{split}
\]
For $1\le p<\infty$ and $0<\alpha\le 1$, the \textit{Besov space} $B_p^\alpha$
is defined as the set of all functions $f\in L^p$ such that
\[
|f|_{B_p^\alpha}:=
\left\{
\begin{array}{ll}
\left(\int_0^{2\pi}\big[s^{-\alpha}\omega_{L^p}^1(f,s)\big]^p\frac{ds}{s}\right)^{1/p}
&
(0<\alpha<1),
\\
\left(\int_0^{2\pi}\big[s^{-1}\omega_{L^p}^2(f,s)\big]^p\frac{ds}{s}\right)^{1/p}
& (\alpha=1)
\end{array}
\right.
\]
is finite. The Besov space is a Banach space under the norm
\[
\|f\|_{B_p^\alpha}:=\|f\|_{L^p}+|f|_{B_p^\alpha}.
\]
These spaces are studied in detail (in a more general setting) in \cite{ST87}
and in many other monographs.
The Riesz projections $P$ and $Q$ are bounded on the Besov spaces $B_p^\alpha$
for $1\le p<\infty$ and $0<\alpha\le 1$ (see \cite[Appendix~2.6]{Pe03}).

Peller proved in the late 1970s that for $1\le p,q<\infty$ and $a\in L^\infty$,
\begin{equation}\label{eq:Peller-0}
Pa\in B_q^{1/q}
\Longleftrightarrow
H(a)\in\cC_q(H^2),
\quad\quad
Qa\in B_p^{1/p}
\Longleftrightarrow
H(\widetilde{a})\in\cC_p(H^2)
\end{equation}
(see \cite[Chap.~6, Theorems~1.1 and 2.1]{Pe03}). From those proofs one can see
that there exist positive constants $c_1$ and $c_2$ depending only on $p$ and $q$
such that
\[
\begin{split}
&
c_1\|Pa\|_{B_q^{1/q}}\le \|H(a)\|_{\cC_q(H^2)}\le c_2\|Pa\|_{B_q^{1/q}},
\\
&
c_1\|Qa\|_{B_p^{1/p}}\le \|H(\widetilde{a})\|_{\cC_p(H^2)}\le c_2\|Qa\|_{B_p^{1/p}}.
\end{split}
\]
From this result and Theorem~\ref{th:Krein-invertibility}(a) we get
the following.
\begin{corollary}\label{co:Krein-representation}
If $1\le p,q<\infty$, then
\[
\begin{split}
& K_{p,0}^{1/p,0}
=
\big\{a\in L^\infty :  Qa\in B_p^{1/p}\big\}
=
L^\infty\cap(B_p^{1/p}+H^\infty),
\\
& K_{0,q}^{0,1/q}
=
\big\{a\in L^\infty : Pa\in B_q^{1/q}\big\}
=
L^\infty\cap(B_q^{1/q}+\overline{H^\infty}),
\\
& K_{p,q}^{1/p,1/q}
=
\big\{a\in L^\infty : Qa\in B_p^{1/p},
Pa\in B_q^{1/q}\big\}
=
L^\infty
\cap
(B_p^{1/p}+H^\infty)
\cap
(B_q^{1/q}+\overline{H^\infty}).
\end{split}
\]
The norms
\[
\|a\|_{L^\infty}+\|Qa\|_{B_p^{1/p}},
\quad
\|a\|_{L^\infty}+\|Pa\|_{B_q^{1/q}},
\quad
\|a\|_{L^\infty}+\|Qa\|_{B_p^{1/p}}+\|Pa\|_{B_q^{1/q}}
\]
are equivalent norms in $K_{p,0}^{1/p,0}$, $K_{0,q}^{0,1/q}$, and
$K_{p,q}^{1/p,1/q}$, respectively.
\end{corollary}
\subsection{Generalized Krein algebras}
We are going to extend the notion of Krein algebras
$K_{p,0}^{1/p,0}$, $K_{0,q}^{0,1/q}$, and $K_{p,q}^{1/p,1/q}$.
Now we take a generalization of the results of
Corollary~\ref{co:Krein-representation} as a definition.

Assume that $1<p,q<\infty$ and $0<\alpha,\beta<1$. Define
\[
\begin{split}
& K_{p,0}^{\alpha,0}
:=
\big\{a\in L^\infty :  Qa\in B_p^\alpha\big\}
=
L^\infty\cap(B_p^\alpha+H^\infty),
\\
& K_{0,q}^{0,\beta}
:=
\big\{a\in L^\infty : Pa\in B_q^\beta\big\}
=
L^\infty\cap(B_q^\beta+\overline{H^\infty}),
\\
& K_{p,q}^{\alpha,\beta}
:=
\big\{a\in L^\infty : Qa\in B_p^\alpha, Pa\in B_q^\beta\big\}
=
L^\infty
\cap
(B_p^\alpha+H^\infty)
\cap
(B_q^\beta+\overline{H^\infty}).
\end{split}
\]
\begin{theorem}[Main result 1]
\label{th:algebra}
Let $1<p,q<\infty$ and $0<\alpha,\beta<1$.
\begin{enumerate}
\item[(a)]
If $\alpha\ge 1/p$, then $K_{p,0}^{\alpha,0}$ is a Banach algebra under the
quasi-submultiplicative norm
\begin{equation}\label{eq:def-gen-Krein1}
\|a\|_{K_{p,0}^{\alpha,0}}:=\|a\|_{L^\infty}+\|Qa\|_{B_p^\alpha}.
\end{equation}

\item[(b)]
If $\beta\ge 1/q$, then $K_{0,q}^{0,\beta}$ is a Banach algebra under the
quasi-submultiplicative norm
\begin{equation}\label{eq:def-gen-Krein2}
\|a\|_{K_{0,q}^{0,\beta}}:=\|a\|_{L^\infty}+\|Pa\|_{B_q^\beta}.
\end{equation}

\item[(c)]
If $\alpha>1/p$, or $\beta>1/q$, or $\alpha=1/p$ and $\beta=1/q$, then
$K_{p,q}^{\alpha,\beta}$ is a Banach algebra under the quasi-submultiplicative norm
\begin{equation}\label{eq:def-gen-Krein3}
\|a\|_{K_{p,q}^{\alpha,\beta}}:=\|a\|_{L^\infty}+\|Qa\|_{B_p^\alpha}+\|Pa\|_{B_q^\beta}.
\end{equation}
\end{enumerate}
\end{theorem}
\begin{theorem}[Main result 2]
\label{th:invertibility}
Let $1<p,q<\infty$, $0<\alpha,\beta<1$,
$1/p+1/q=\alpha+\beta\in(0,1]$.
\begin{enumerate}
\item[(a)]
Suppose $\alpha\ge 1/p$ and $K$ is either $K_{p,0}^{\alpha,0}$ or
$K_{p,q}^{\alpha,\beta}$. If $a\in K$, then
\[
a\in GK \Longleftrightarrow a\in G(C+H^\infty).
\]

\item[(b)]
Suppose $\beta\ge 1/q$ and $K$ is either $K_{0,q}^{0,\beta}$ or
$K_{p,q}^{\alpha,\beta}$. If $a\in K$, then
\[
a\in GK \Longleftrightarrow a\in G(C+\overline{H^\infty}).
\]
\end{enumerate}
\end{theorem}
These statements were proved in \cite[Chap.~4]{BS83} for the particular
case of $K_{p,q}^{\alpha,\beta}$ in which the parameters satisfy
$1/p+1/q=\alpha+\beta=1$. Those proofs are based on Krein's ideas
\cite{Kr66}. Since the book \cite{BS83} is no longer available to a wide audience,
we decided to present self-contained proofs of these results.
\subsection{Szeg\H{o}-Widom type limit theorem}
Let $N$ be a natural number. For a Banach space $X$, let $X_N$ and
$X_{N\times N}$ be the spaces of vectors and matrices with entries in $X$.
The operators $I$, $J$, $P$, and $Q$ are defined on vector spaces elementwise;
the Hankel operators on $H_N^2$ are defined for
$a\in L_{N\times N}^\infty$ in the same way as in (\ref{eq:def-Hankel})
and the Toeplitz operators are defined by
\[
T(a):H_N^2\to H_N^2,
\quad
f\mapsto PM(a)Pf,
\quad
T(\widetilde{a}):H_N^2\to H_N^2,
\quad
f\mapsto JQM(a)QJf.
\]
The matrix of the Toeplitz operator $T(a)$ in the standard basis of the space $H_N^2$
is the infinite Toeplitz matrix
\[
\left(
\begin{array}{cccc}
a_0 & a_{-1} & a_{-2} & \dots\\
a_1 & a_0 & a_{-1} & \dots\\
a_2 & a_1 & a_0 & \dots\\
\dots & \dots & \dots & \dots
\end{array}
\right)
\]
where $\{a_k\}_{k\in\Z}$ is the sequence of $N\times N$ matrices which are the
Fourier coefficients of the generating function (symbol) $a\in L_{N\times N}^\infty$.
Let $T_n(a)=(a_{j-k})_{j,k=0}^n$. This $nN \times nN$ matrix
is usually called a block Toeplitz matrix. We denote its determinant by $D_n(a)$.

In 1915, Szeg\H{o} proved that $D_{n-1}(a)/D_n(a)$ tends to the geometric mean
$G(a)$ of $a$ if $a$ is a scalar nonnegative function such that $a\in L^1$ and
$\log a\in L^1$. This result, now called the first Szeg\H{o} limit theorem,
has been subsequently extended into different directions. We will not
go into details, but notice that Krein \cite{Kr66} observed that his algebra
$K_{2,2}^{1/2,1/2}$ can be useful in asymptotic analysis of Toeplitz determinants.
The most general results for the case where positivity is replaced by some kind of
sectoriality are those of Krein and Spitkovsky \cite{KS83}.

In 1952, Szeg\H{o} proved his second (strong) limit theorem, which says that
if $a$ is a positive scalar function with H\"older continuous derivative, then
$D_n(a)\sim G(a)^{n+1}E(a)$, where $E(a)$ is some completely determined nonzero
constant. It was Widom \cite{W76} who proved the strong Szeg\H{o} theorem in the
block case $N>1$ for the first time under the assumption that
$a\in (K_{2,2}^{1/2,1/2})_{N\times N}$.

We suppose that the reader is familiar with basic facts on Fredholm operators
and with properties of Schatten-von Neumann classes and (regularized) operator
determinants (otherwise consult Section~\ref{sec:preliminaries}).

Functions in $K_{p,q}^{\alpha,\beta}$ may be discontinuous, hence some care is
needed in the definition of the ``geometric mean" $G(a)$.
For $a\in L_{N\times N}^\infty$, we denote by $h_ra$ the harmonic extension,
\[
h_r a(e^{i\theta}):=\sum_{n=-\infty}^{\infty}a_nr^{|n|}e^{in\theta}
\quad (0\le r<1,\quad 0\le\theta<2\pi).
\]
\begin{lemma}\label{le:constant-G}
Let $a\in (C+H^\infty)_{N\times N}$ or $a\in (C+\overline{H^\infty})_{N\times N}$.
If $T(a)$ is Fredholm of index zero on $H_N^2$, then the limit
\begin{equation}\label{eq:def-G}
G(a):=
\lim_{r\to 1-0}\exp\left(\frac{1}{2\pi}\int_0^{2\pi}
\log\det h_ra(e^{i\theta})d\theta\right)
\end{equation}
exists, is finite and nonzero.
\end{lemma}
The proof of this lemma is given in \cite[Proposition~10.6(a)]{BS06}
for $(C+H^\infty)_{N\times N}$ and it works equally also for
$(C+\overline{H^\infty})_{N\times N}$.
\begin{theorem}\label{th:Szego-Widom}
Let $K$ be one of the algebras $K_{1,0}^{1,0}$, $K_{0,1}^{0,1}$, or
$K_{p,q}^{\alpha,\beta}$ with $p,q,\alpha,\beta$ satisfying
\[
1<p,q<\infty,
\quad
0<\alpha,\beta<1,
\quad
1/p+1/q=\alpha+\beta=1,
\quad
-1/2<\alpha-1/p<1/2.
\]
If $a\in K_{N\times N}$
and $T(a)$ is Fredholm of index zero on $H_N^2$, then $a^{-1}\in K_{N\times N}$,
the operator
\[
H(a)H(\widetilde{a}^{-1})=I-T(a)T(a^{-1})
\]
is of trace class on $H_N^2$, and
\[
\lim_{n\to\infty}\frac{D_n(a)}{G(a)^{n+1}}=\det T(a)T(a^{-1}),
\]
where the last $\det$ refers to the determinant defined for operators
differing from the identity by an operator of trace class.
\end{theorem}
For $1/p=1/q=\alpha=\beta=1/2$ this result goes back to Widom \cite{W76}.
Two of the authors proved it in this form in \cite[Theorem~6.14]{BS83}.
The proof is also reproduced in \cite[Theorem~10.32]{BS06}.

Alternative proofs of the strong Szeg\H{o} and Szeg\H{o}-Widom limit theorems
were found later. We refer to \cite{BG05,BS83,BS99,BS06,BW06,E03,GS58,S05,W76}
and the references given there for more complete information.
\subsection{Higher order asymptotic formulas for block Toeplitz determinants}
For $a\in L_{N\times N}^\infty$ and $n\in\Z_+:=\{0,1,2,\dots\}$
define the operators $P_n$ and $Q_n$ on $H_N^2$ by
\[
P_n:\sum_{k=0}^\infty a_k\chi_k\mapsto \sum_{k=0}^na_k\chi_k,
\quad
Q_n:=I-P_n.
\]
The operator $P_nT(a)P_n:P_nH_N^2\to P_nH_N^2$ may be identified with the finite
block Toeplitz matrix $T_n(a):=(a_{j-k})_{j,k=0}^n$. Let $W$ be the Wiener algebra
of functions $a:\T\to\C$ with absolutely convergent Fourier series. For generalized
Krein algebras, define the ``conjugation number" by
\begin{equation}\label{eq:def-lambda}
\lambda:=\left\{\begin{array}{lllll}
1/p & \mbox{for} & K_{p,0}^{\alpha,0} &\mbox{with}& \alpha\ge 1/p,
\\[2mm]
1/q & \mbox{for} & K_{0,q}^{0,\beta}  &\mbox{with}& \beta\ge 1/q,
\\[2mm]
1/p+1/q &\mbox{for} & K_{p,q}^{\alpha,\beta}&\mbox{with}& 1/p+1/q=\alpha+\beta.
\end{array}\right.
\end{equation}
\begin{theorem}[Main result 3]
\label{th:main}
Let $1<p,q<\infty$ and $0<\alpha,\beta<1$. Suppose $K$ is one of the algebras
\[
W\cap K_{p,0}^{1/p,0},
\quad
K_{p,0}^{\alpha,0} \mbox{ with } \alpha>1/p,
\quad
W\cap K_{0,q}^{0,1/q},
\quad
K_{0,q}^{0,\beta} \mbox{ with } \beta>1/q,
\quad
W\cap K_{p,q}^{1/p,1/q},
\]
or $K_{p,q}^{\alpha,\beta}$ with $\alpha\ne 1/p$, $1/p+1/q=\alpha+\beta\in(0,1)$,
$-1/2<\alpha-1/p<1/2$. If $a\in K_{N\times N}$ and both $T(a)$ and
$T(\widetilde{a})$ are invertible on $H_N^2$, then the following statements hold.
\begin{enumerate}
\item[(a)]
The matrix function $a$ is invertible in $K_{N\times N}$ and
admits canonical right and left Wiener-Hopf
factorizations in $K_{N\times N}$, that is,
there exist
\[
u_-,v_-\in G(K\cap \overline{H^\infty})_{N\times N},
\quad
u_+,v_+\in G(K\cap H^\infty)_{N\times N}
\]
such that $a=u_-u_+=v_+v_-$.

\item[(b)]
Let $\lambda$ be defined by {\rm(\ref{eq:def-lambda})} and $m$ be the smallest
integer such that $1\le \lambda m$. If
\begin{equation}\label{eq:def-bc}
b:=v_-u_+^{-1},
\quad
c:=u_-^{-1}v_+,
\end{equation}
where $u_\pm,v_\pm$ are defined in part {\rm(a)},
then the operators
$H(\widetilde{c})H(b)$ and $H(b)H(\widetilde{c})$ belong to the
Schatten-von Neumann class $\cC_m(H_N^2)$.

\item[(c)] Let $b,c$, and $m$ be defined as in part {\rm(b)}. Then
\begin{equation}\label{eq:formula}
\lim_{n\to\infty}\frac{D_n(a)}{G(a)^{n+1}}
\exp\left\{
-\sum_{j=1}^{m-1}\frac{1}{j}\,\tr\!\left[
\left(
\sum_{k=0}^{m-1}F_{n,k}
\right)^j
\right]
\right\}
=
\frac{1}{\det_m T(\widetilde{c})T(\widetilde{b})},
\end{equation}
where
\begin{equation}\label{eq:def-F}
F_{n,k} :=
P_nT(c)Q_n\big(Q_nH(b)H(\widetilde{c})Q_n\big)^kQ_nT(b)P_n
\quad
(k\in\Z_+)
\end{equation}
and $\det_m$ denotes the $m$-regularized operator determinant defined for
operators differing from the identity by an operator in $\cC_m(H_N^2)$
{\rm(}see Section~{\rm\ref{sec:determinant})}.
\end{enumerate}
\end{theorem}
Under the assumptions of Theorem~\ref{th:main} we cannot guarantee that
$T(a)T(a^{-1})-I$ is of trace class, which implies that the arguments of the
proof of Theorem~\ref{th:Szego-Widom} do not work in this case. However,
we can guarantee that $H(\widetilde{c})H(b)$ and $H(b)H(\widetilde{c})$
belong to $\cC_m(H_N^2)$ where $m>1$. In this case the asymptotic formulas
for block Toeplitz determinants are subject to higher order corrections
involving additional terms and regularized operator determinants.

Notice also that similar results were obtained in \cite[Section 4]{BS80},
\cite[Theorem~6.20]{BS83}, \cite[Theorem~10.37]{BS06}, \cite[Theorem~20]{K06}
for weighted Wiener algebras and H\"older spaces $C^\gamma$. For Wiener
algebras with power weights and for $C^\gamma$,  the formula in (c) is a
little bit simpler, because it does not contain $F_{n,m-1}$. In contrast to the case
of weighted Wiener algebras and $C^\gamma$, which consist of continuous
functions only, the generalized Krein algebras $K_{p,q}^{\alpha,\beta}$ may
contain discontinuous functions.
\subsection{About this paper}
Section~\ref{sec:preliminaries} contains operator-theoretic preliminaries.
In Section~\ref{sec:Krein}, we slightly extend Krein's results \cite{Kr66}
on Banach algebras generated by ideals. Section~\ref{sec:aux} contains
auxiliary results from the theory of Besov spaces, algebras of multiplication
operators on weighted $\ell_2$-spaces, as well as Peller's results
on the boundedness, compactness, and Schatten-von Neumann behavior of Hankel-type
operators. Let $c^\gamma$ be the closure of the set of all Laurent polynomials
in the norm of the H\"older space $C^\gamma$, $0<\gamma<1$.
In Section~\ref{sec:CH-and-cH}, we show that $C^\gamma+H^\infty$ and
$c^\gamma+H^\infty$ are Banach algebras and that $c^\gamma+H^\infty$ is inverse
closed in $C+H^\infty$. Section~\ref{sec:gen-Krein-algebras} contains the
proofs of Theorems~\ref{th:algebra} and \ref{th:invertibility}. The proof
of Theorem~\ref{th:invertibility} uses essentially the inverse closedness of
$c^\gamma+H^\infty$ in $C+H^\infty$.
In Section~\ref{sec:proof-asymptotic},
we prove Theorem~\ref{th:main} by first applying the factorization theory developed
in \cite{HS84} to generalized Krein algebra and then using an abstract
higher order asymptotic formula for Toeplitz determinants \cite[Theorem~15]{K06},
which is contained implicitly in \cite[Theorem~6.20]{BS83} and
\cite[Theorem~10.37]{BS06}.
\section{Operator-theoretic preliminaries}\label{sec:preliminaries}
\subsection{Commutative Banach algebras}
Let $\cA$ be an algebra. A subalgebra $\cJ$ of $\cA$ is called an algebraic
\textit{two-sided ideal} of $\cA$ if $aj\in\cJ$ and $ja\in\cJ$ for all
$a\in\cA$ and $j\in\cJ$.
Given two Banach algebras $\cA$ and $\cB$, a map $\varphi:\cA\to\cB$ is called
a Banach algebra homomorphism if $\varphi$ is a bounded linear operator and
$\varphi(ab)=\varphi(a)\varphi(b)$ for all $a,b\in\cA$.
Now let $\cA$ be a commutative Banach algebra with identity element  $e$.
The Banach algebra homomorphisms of $\cA$ into $\C$ which send $e$ to $1$ are
called multiplicative linear functionals of $\cA$. A proper closed two-sided
ideal $\cJ$ of $\cA$ is called a \textit{maximal ideal} if it is not properly
contained in any other proper closed two-sided ideal of $\cA$. Let $\cM_\cA$ denote the
set of all maximal ideals of $\cA$ and let $M_\cA$ stand for the set of all
multiplicative linear functionals of $\cA$. One can show that the map
$M_\cA\to\cM_\cA$, $\varphi\to\Ker\varphi$ is bijective. Therefore no distinction
is usually made between multiplicative linear functionals and maximal ideals.

The formula $\widehat{a}(m)=m(a)$ ($m\in M_\cA$) assigns a function
$\widehat{a}:M_\cA\to\C$ to each $a\in\cA$. Let $\widehat{A}$ be the set
$\{\widehat{a}:a\in\cA$\}. The Gelfand topology on $M_\cA$ is the coarsest
topology on $M_\cA$ which makes all functions $\widehat{a}\in\widehat{A}$
continuous. The set $M_\cA$ equipped with the Gelfand topology is called the
maximal ideal space of $\cA$.
\begin{theorem}[Gelfand]
Let $\cA$ be a commutative Banach algebra with identity element and let
$M_\cA$ be the maximal ideal space of $\cA$. An element $a\in\cA$
is invertible if and only if $\widehat{a}(m)\ne 0$ for all $m\in M_\cA$.
\end{theorem}
A proof of this result is in every textbook on Banach algebras.
\subsection{Fredholm operators}
The facts collected in this subsection can be found in most textbooks on operator
theory (for instance, in  \cite{D98,GGK90,GK57}).
Let $X$ be a Banach space, $\cB(X)$ be the Banach algebra of all bounded
linear operators on $X$, $\cC_0(X)$ be the set of all finite-rank operators,
and $\cC_\infty(X)$ be the closed two-sided ideal of all compact operators on $X$.
For $A\in\cB(X)$, put $\Ker A=\{x\in X: Ax=0\}$ and $\im A=AX$. The operator
$A$ is said to be \textit{Fredholm} (on $X$) if $\im A$ is closed in $X$ and both
$\dim\Ker A$ and $\dim (X/\im A)$ are finite. The integer
\[
\ind A:=\dim\Ker A-\dim(X/\im A)
\]
is then referred to as the \textit{index} of $A$.

If $A$ and $B$ are Fredholm, then $AB$ is also Fredholm and $\ind (AB)=\ind A+\ind B$.
If $A$ is Fredholm and $K\in\cC_\infty(X)$, then $A+K$ is Fredholm and
$\ind (A+K)=\ind A$. An operator $A\in\cB(X)$ is Fredholm if and only if there
exists an operator $R\in\cB(X)$ such that $AR-I\in\cC_\infty(X)$ and
$RA-I\in\cC_\infty(X)$. Such an operator $R$ is called a \textit{regularizer}
of the operator $A$.
\subsection{Schatten-von Neumann ideals}\label{sec:Schatten}
All facts stated in the rest of this section are proved in \cite[Chap.~3--4]{GK69}.
Let $H$ be a separable Hilbert space. Given an operator $A\in\cB(H)$ define for
$n\in\Z_+$,
\[
s_n(A):=\inf\{\|A-F\|_{\cB(H)}\ :\ F\in\cC_0(H),\ \dim F(H)\le n\}.
\]
For $1\le p<\infty$, the collection of all operators $K\in\cB(H)$ satisfying
\begin{equation}\label{eq:Schatten-norm}
\|K\|_{\cC_p(H)}:=\Big(\sum_{n\in\Z_+}s_n^p(K)\Big)^{1/p}<\infty
\end{equation}
is denoted by $\cC_p(H)$ and referred to as a \textit{Schatten-von Neumann class}.
This is a Banach space under the norm (\ref{eq:Schatten-norm}).
Note that $\cC_\infty(H)=\{K\in\cB(H):s_n(K)\to 0\mbox{ as }n\to\infty\}$ and
\[
\|K\|_{\cC_\infty(H)}=\sup_{n\in\Z_+}s_n(K)=\|K\|_{\cB(H)}.
\]
The operators belonging to $\cC_1(H)$ are called \textit{trace class operators}.
\begin{lemma}\label{le:Schatten-properties}
Let $1\le p,q,r\le\infty$.
\begin{enumerate}
\item[(a)]
If $p<q$ and $K\in\cC_p(H)$, then $K\in\cC_q(H)$ and $\|K\|_{\cC_q(H)}\le\|K\|_{\cC_p(H)}$.
\item[(b)]
If $A\in\cB(H)$ and $K\in\cC_p(H)$, then $AB,BA\in\cC_p(H)$ and
\[
\max\{\|AK\|_{\cC_p(H)},\|KA\|_{\cC_p(H)}\}
\le
\|K\|_{\cC_p(H)}\|A\|_{\cB(H)}.
\]
\item[(c)]
If $1/r=1/p+1/q$ and $K\in\cC_p(H)$, $L\in\cC_q(H)$, then $KL\in\cC_r(H)$ and
\[
\|KL\|_{\cC_r(H)}\le\|K\|_{\cC_p(H)}\|L\|_{\cC_q(H)}.
\]
\end{enumerate}
\end{lemma}
Lemma~\ref{le:Schatten-properties}(b) implies that $\cC_p(H)$ is a two-sided
ideal of $\cB(H)$. This ideal is not closed in $\cB(H)$ if $p<\infty$ and ${\rm dim}\,H=\infty$.
\subsection{Regularized operator determinants}\label{sec:determinant}
Let $A\in\cB(H)$ be an operator of the form $I+K$ with $K\in\cC_1(H)$. If
$\{\lambda_j(K)\}_{j\ge 0}$ denotes the sequence of the nonzero eigenvalues of
$K$ counted up to algebraic multiplicity, then the product
$\prod_{j\ge 0}(1+\lambda_j(K))$ is absolutely convergent. The \textit{determinant}
of $A$ is defined by
\[
\det A=\det(I+K)=\prod_{j\ge 0}(1+\lambda_j(K)).
\]
If $K\in\cC_m(H)$, where $m>1$ is an integer, one can still define a determinant
of $I+K$, but for ideals larger than $\cC_1(H)$, the above definition requires a
regularization. A simple computation (see \cite[Lemma~6.1]{S77}) shows that then
\[
R_m(K):=(I+K)\exp\Bigg(\sum_{j=1}^{m-1}\frac{(-K)^j}{j}\Bigg)-I\in\cC_1(H).
\]
Thus, it is natural to define $\det_1(I+K):=\det(I+K)$ and
$\det_m(I+K):=\det(I+R_m(K))$ for $m>1$. One calls
$\det_m(I+K)$ the $m$-\textit{regularized determinant} of $A=I+K$.

Regularized operator determinants have some useful properties. For instance,
the operator $I+K$ is invertible on $H$ if and only if $\det_m(I+K)\ne 0$.
\section{Banach algebras generated by ideals}\label{sec:Krein}
\subsection{Subalgebras generated by ideals}
The results of this section are essentially due to Krein \cite{Kr66}. They
appeared in this form in \cite[Chap.~4]{BS83}. Since both sources may not be
available to a wide audience, it seems reasonable to give here complete proofs.
\begin{lemma}\label{le:subalgebras}
Let $\cA$ be an algebra and $e\in\cA$ be the identity element. Suppose
$\cJ_1\subset\cA$ and $\cJ_2\subset\cA$ are two-sided ideals in $\cA$ and
$\cL$ is a subalgebra of $\cA$. Let
\begin{equation}\label{eq:subalgebras-1}
p\in\cA,\quad p^2=p,\quad q:=e-p.
\end{equation}
Then the sets
\begin{equation}\label{eq:subalgebras-2}
\cL_1:=\{a\in\cL\ :\ paq\in\cJ_1\},
\quad
\cL_2:=\{a\in\cL\ :\ qap\in\cJ_2\},
\quad
\cL_*:=\cL_1\cap\cL_2
\end{equation}
are subalgebras of $\cL$.
\end{lemma}
\begin{proof}
Because $peq=p(e-p)=p-p^2=0$ and $qep=(e-p)p=p-p^2=0$,
we have $e\in\cL_1$, $e\in\cL_2$, and thus $e\in\cL_1\cap\cL_2$.
If $a,b\in\cL_1$, then $paq\in\cJ_1$ and $pbq\in\cJ_1$. We have
\[
pabq=pa(p+q)bq=pap\cdot pbq+paq\cdot qbp.
\]
Since $\cJ_1$ is a two-sided ideal of $\cL$ and $pap\in\cL$, $qbp\in\cL$,
we conclude that $pabq\in\cJ_1$, that is, $ab\in\cL_1$. Thus, $\cL_1$
is an algebra with identity $e\in\cL$. Analogously one can prove that
$\cL_2$ is a subalgebra of $\cL$. The statement for $\cL_*=\cL_1\cap\cL_2$
is now trivial.
\end{proof}
\subsection{Banach algebras generated by complete normed ideals}
Throughout this section we assume that $\cA$ is a Banach algebra with identity
$e$ and a Banach algebra norm $\|\cdot\|$.
\begin{theorem}\label{th:Krein1}
Let $\cJ_1\subset\cA$, $\cJ_2\subset\cA$ be two-sided ideals of $\cA$, which
become Banach spaces under norms $\|\cdot\|_{\cJ_i}$, $i=1,2$. Assume that for
every $i=1,2$ and every $x\in \cJ_i$,
\begin{equation}\label{eq:Krein1-1}
\|x\|\le \|x\|_{\cJ_i}.
\end{equation}
If $\cL$ is a closed subalgebra of $\cA$, then $\cL_1$, $\cL_2$, and $\cL_*$
defined by {\rm (\ref{eq:subalgebras-1})--(\ref{eq:subalgebras-2})} are
Banach spaces under the norms
\[
\begin{split}
\|a\|_1 &:=\big(\|pap\|^2+\|qaq\|^2+\|paq\|_{\cJ_1}^2+\|qap\|^2\big)^{1/2},
\\
\|a\|_2 &:=\big(\|pap\|^2+\|qaq\|^2+\|paq\|^2+\|qap\|_{\cJ_2}^2\big)^{1/2},
\\
\|a\|_* &:=\big(\|pap\|^2+\|qaq\|^2+\|paq\|_{\cJ_1}^2+\|qap\|_{\cJ_2}^2\big)^{1/2},
\end{split}
\]
respectively.
\end{theorem}
\begin{proof}
It is easy to see that $\|\cdot\|_i$ and $\|\cdot\|_*$ are norms.
The triangle inequality for $\|\cdot\|_i$ as well as for $\|\cdot\|_*$ follows
from the triangle inequality for $\|\cdot\|$ and $\|\cdot\|_{\cJ_i}$ and from
the Minkowski inequality for sums.

If $a\in\cL_1$, then $paq\in\cJ_1$. In view of (\ref{eq:Krein1-1}), we have
\begin{equation}\label{eq:Krein1-3}
\begin{split}
\|a\|
& \le
\|pap\|+\|qaq\|+\|paq\|+\|qap\|
\\
&\le
\|pap\|+\|qaq\|+\|paq\|_{\cJ_1}+\|qap\|
\\
&\le
4\big(\|pap\|^2+\|qaq\|^2+\|paq\|_{\cJ_1}^2+\|qap\|^2\big)^{1/2}
\\
& =4\|a\|_1.
\end{split}
\end{equation}
It is obvious that
\begin{equation}\label{eq:Krein1-4}
\|paq\|_{\cJ_1}\le\|a\|_1.
\end{equation}
Let $\{a_n\}\subset\cJ_1$ be a Cauchy sequence in $\cL_1$, that is, $\|a_n-a_m\|_1\to 0$
as $n,m\to\infty$. From (\ref{eq:Krein1-3}) it follows that $\{a_n\}$ is a
Cauchy sequence in $\cL$. Since $\cL$ is closed, there exists an $a\in\cL$
such that $\|a-a_n\|\to 0$ as $n\to\infty$. This implies that as $n\to\infty$,
\begin{eqnarray}
\label{eq:Krein1-5}
\|pap-pa_np\| &\to& 0,
\\
\label{eq:Krein1-6}
\|paq-pa_nq\| &\to& 0,
\\
\label{eq:Krein1-7}
\|qap-qa_np\| &\to& 0,
\\
\label{eq:Krein1-8}
\|qaq-qa_nq\| &\to& 0.
\end{eqnarray}
On the other hand, from (\ref{eq:Krein1-4}) we see that
$\|pa_nq-pa_mq\|_{\cJ_1}\le\|a_n-a_m\|_1\to 0$ as $n,m\to\infty$,
that is, $\{pa_nq\}$ is a Cauchy sequence in $\cJ_1$. By the hypothesis,
$(\cJ_1,\|\cdot\|_{\cJ_1})$ is a Banach space. Thus there is an element
$c\in\cJ_1$ such that $\|c-pa_nq\|_{\cJ_1}\to 0$ as $n\to\infty$. In view
of (\ref{eq:Krein1-1}), this gives
\begin{equation}\label{eq:Krein1-9}
\|c-pa_nq\|\to 0\quad\mbox{as}\quad n\to\infty.
\end{equation}
From (\ref{eq:Krein1-6}) and (\ref{eq:Krein1-9}) it follows that $c=paq$.
Hence $a\in\cL_1$ and
\begin{equation}\label{eq:Krein1-10}
\|paq-pa_nq\|_{\cJ_1}\to 0\quad\mbox{as}\quad n\to\infty.
\end{equation}
Combining (\ref{eq:Krein1-5}), (\ref{eq:Krein1-7}), (\ref{eq:Krein1-8}),
and (\ref{eq:Krein1-10}), we get $\|a-a_n\|_1\to 0$ as $n\to\infty$.
Thus $a\in\cL_1$ and $\cL_1$ is closed.

It is clear that the same argument works also for $\cL_2$ and $\cL_*$.
Thus $\cL_2$ and $\cL_*$ are closed, too.
\end{proof}
\begin{theorem}\label{th:Krein2}
Under the assumptions of Theorem~{\rm\ref{th:Krein1}}, the norms
\[
\begin{split}
\|a\|_1' &:=\|a\|+\|paq\|_{\cJ_1}\quad (a\in\cL_1),
\\
\|a\|_2' &:=\|a\|+\|qap\|_{\cJ_2}\quad (a\in\cL_2),
\\
\|a\|_*' &:=\|a\|+\|paq\|_{\cJ_1}+\|qap\|_{\cJ_2}\quad (a\in\cL_*)
\end{split}
\]
are equivalent to the norms $\|a\|_1$, $\|a\|_2$, and $\|a\|_*$, respectively.

If, in addition,
\[
\|p\|=\|q\|=1
\]
and for every $i=1,2$, for every $x\in\cJ_i$, and every $a\in\cA$,
\begin{equation}\label{eq:Krein1-2}
\max\{\|ax\|_{\cJ_i},\|xa\|_{\cJ_i}\}\le\|x\|_{\cJ_i}\|a\|,
\end{equation}
then $\|\cdot\|_1'$, $\|\cdot\|_2'$, and $\|\cdot\|_*'$ are Banach algebra norms.
\end{theorem}
\begin{proof}
Let us prove the statement for $\|\cdot\|_*'$. By analogy with (\ref{eq:Krein1-3}),
we get
\begin{equation}\label{eq:Krein2-1}
\|a\|\le 4\|a\|_*.
\end{equation}
It is obvious that
\begin{equation}\label{eq:Krein2-2}
\|paq\|_{\cJ_1}\le\|a\|_*,
\quad
\|qap\|_{\cJ_2}\le\|a\|_*.
\end{equation}
From (\ref{eq:Krein2-1}) and (\ref{eq:Krein2-2}) we see that $\|a\|_*'\le 6\|a\|_*$.

On the other hand,
\[
\begin{split}
\|a\|_*
&\le
\|pap\|+\|qaq\|+\|paq\|_{\cJ_1}+\|qap\|_{\cJ_2}
\\
&\le
(\|p\|^2+\|q\|^2+1)(\|a\|+\|paq\|_{\cJ_1}+\|qap\|_{\cJ_2})
\\
&=
(\|p\|^2+\|q\|^2+1)\|a\|_*',
\end{split}
\]
that is, the norms $\|\cdot\|_*$ and $\|\cdot\|_*'$ are equivalent.

It is clear that $\|e\|_*'=\|e\|+\|pq\|_{\cJ_1}+\|qp\|_{\cJ_2}=\|e\|=1$.
If $\|p\|=\|q\|=1$, then from (\ref{eq:Krein1-2}) it follows that
\begin{equation}\label{eq:Krein2-3}
\begin{split}
\|pabq\|_{\cJ_1}
&=
\|pa(p+q)bq\|_{\cJ_1}
\\
&\le
\|pap\|\,\|pbq\|_{\cJ_1}+\|paq\|_{\cJ_1}\|qbq\|
\\
& \le
\|a\|\,\|pbq\|_{\cJ_1}+\|paq\|_{\cJ_1}\|b\|
\end{split}
\end{equation}
and
\begin{equation}\label{eq:Krein2-4}
\|qabp\|_{\cJ_2}
\le
\|qap\|_{\cJ_2}\|b\|+\|a\|\,\|qbp\|_{\cJ_2}.
\end{equation}
Combining $\|ab\|\le\|a\|\,\|b\|$ and (\ref{eq:Krein2-3})--(\ref{eq:Krein2-4}),
we get
\[
\begin{split}
\|ab\|_*'
&\le
\|a\|\,\|b\|+
\|a\|\,\|pbq\|_{\cJ_1}+\|paq\|_{\cJ_1}\|b\|+
\|qap\|_{\cJ_2}\|b\|+\|a\|\,\|qbp\|_{\cJ_2}
\\
&=
\|a\|\,\|b\|_*'+\|b\|(\|paq\|_{\cJ_1}+\|qap\|_{\cJ_2})
\\
&\le
\|a\|\,\|b\|_*'+\|b\|_*'(\|paq\|_{\cJ_1}+\|qap\|_{\cJ_2})
\\
&=
\|a\|_*'\|b\|_*'
\end{split}
\]
for all $a,b\in\cL_*$. The proof for $\cL_1$ and $\cL_2$ is similar.
\end{proof}
\begin{remark}
One can show that under the assumptions of Theorem~\ref{th:Krein1} and
(\ref{eq:Krein1-2}) one has  $\|ab\|_i\le\|a\|_i\|b\|_i$ for all $a,b\in\cL_i$
and $\|ab\|_*\le\|a\|_*\|b\|_*$ for all $a,b\in\cL_*$. However,
$\|e\|_i=\sqrt{2}$ and $\|e\|_*=\sqrt{2}$ if $\|p\|=\|q\|=1$.
\end{remark}
\subsection{Operator algebras generated by ideals of compact operators}
Let $X$ be a Banach space and $\cP,\cQ\in\cB(X)$ be two complementary projections,
that is, $\cP^2=\cP$ and $\cQ=I-\cP$.
\begin{theorem}\label{th:Krein3}
Let $\cJ_1$, $\cJ_2$ be (not necessarily closed) two-sided ideals in $\cB(X)$
such that $\cC_0(X)\subset\cJ_1\subset\cC_\infty(X)$ and
$\cC_0(X)\subset\cJ_2\subset\cC_\infty(X)$,
let $\cL$ be a (not necessarily closed) subalgebra of $\cB(X)$, and let
\[
\cL_1:=\{A\in\cL: \cP A\cQ\in\cJ_1\},
\quad
\cL_2:=\{A\in\cL: \cQ A\cP\in\cJ_2\},
\quad
\cL_*:=\cL_1\cap\cL_2.
\]
\begin{enumerate}
\item[(a)]
If $A\in\cL_1$ is invertible and $\cP A\cP|\im\cP$ is Fredholm on $\cP X$,
then $A^{-1}\in\cL_1$.

\item[(b)]
If $A\in\cL_2$ is invertible and $\cQ A\cQ|\im\cQ$ is Fredholm on $\cQ X$,
then $A^{-1}\in\cL_2$.

\item[(c)]
If $A\in\cL_*$ is invertible, then $A^{-1}\in\cL_*$.
\end{enumerate}
\end{theorem}
\begin{proof}
(a) Put $A_1:=\cP A\cP+\cQ A\cQ$. From the invertibility of $A$ and
\[
A-A_1-\cQ A\cP=\cP A\cQ\in\cJ_1\subset\cC_\infty(X)
\]
it follows that $A_1+\cQ A\cP$ is Fredholm and has index zero.

Let $C\in\cB(\cP X)$ be a regularizer of the operator $D:=\cP A\cP|\im\cP$, which is
Fredholm on $\cP X$. Then $CD=\cP+\cP T\cP$, where $T\in\cC_\infty(\cP X)$.
In that case
\begin{equation}\label{eq:Krein3-1}
(I+\cQ A\cP C\cP)(\cP A\cP+\cQ A\cQ)
=
A_1+\cQ A\cP CD
=
A_1+\cQ A\cP+T_1,
\end{equation}
where $T_1\in\cC_\infty(X)$. The operator $A_1+\cQ A\cP+T_1$ is Fredholm
and has index zero because $A_1+\cQ A\cP$ is so. On the other hand, it is
easy to check that
\begin{equation}\label{eq:Krein3-2}
(I+\cQ A\cP C\cP)^{-1}=I-\cQ A\cP C\cP.
\end{equation}
Thus from (\ref{eq:Krein3-1}) and (\ref{eq:Krein3-2}) it follows that the
operator $A_1$ is Fredholm and has index zero. Consequently, there exist
closed subspaces $M,N\subset X$ such that
\[
X=\Ker A_1\oplus M,\quad X=\im A_1\oplus N.
\]
Hence $\dim\Ker A_1=\dim N$.

Denote by $\widehat{A_1}$ the restriction of $A$ to $M$. Then $\widehat{A_1}:M\to\im A_1$
is continuous, one-to-one, and hence its continuous inverse $(\widehat{A_1})^{-1}$ exists.
Define $B_1\in\cB(X)$ by
\[
B_1x=\left\{
\begin{array}{lll}
(\widehat{A_1})^{-1}x &\mbox{if}& x\in\im A_1,\\
0 &\mbox{if} & x\in N.
\end{array}\right.
\]
Denote the projection onto $\Ker A_1$ parallel to $M$ by $Q_1$. We have for every
$x\in X$,
\[
B_1A_1x=B_1A_1(x-Q_1x)=B_1\widehat{A_1}(x-Q_1x)=x-Q_1x,
\]
that is,
\begin{equation}\label{eq:Krein3-3}
B_1A_1=I-Q_1.
\end{equation}
Set
\begin{equation}\label{eq:Krein3-4}
\begin{array}{llll}
A_{11}:=\cP A\cP,
&
A_{12}:=\cP A\cQ,
&
A_{21}:=\cQ A\cP,
&
A_{22}:=\cQ A\cQ,
\\[2mm]
B_{11}:=\cP A^{-1}\cP,
&
B_{12}:=\cP A^{-1}\cQ,
&
B_{21}:=\cQ A^{-1}\cP,
&
B_{22}:=\cQ A^{-1}\cQ.
\end{array}
\end{equation}
From $AA^{-1}=I$ we get
\begin{equation}\label{eq:Krein3-5}
0=\cP\cQ=\cP AA^{-1}\cQ=\cP A(\cP+\cQ)A^{-1}\cQ=A_{11}B_{12}+A_{12}B_{22}.
\end{equation}
It is obvious that
\begin{equation}\label{eq:Krein3-6}
A_1=A_{11}+A_{12},\quad A_{22}B_{12}=0.
\end{equation}
Combining (\ref{eq:Krein3-5}) and (\ref{eq:Krein3-6}), we get
$A_1B_{12}+A_{12}B_{22}=0$. Multiplying this equality by $B_1$ from the left
and taking into account (\ref{eq:Krein3-3}), we get
$(I-Q_1)B_{12}+B_1A_{12}B_{22}=0$. Hence
\begin{equation}\label{eq:Krein3-7}
B_{12}=Q_1B_{12}-B_1A_{12}B_{22}.
\end{equation}
By the hypothesis, $A\in\cL_1$, whence $A_{12}=\cP A\cQ\in\cJ_1$. Thus,
\begin{equation}\label{eq:Krein3-8}
B_1A_{12}B_{22}\in\cJ_1.
\end{equation}
Since $A_1$ is Fredholm, the projection $Q_1$ has finite rank. Hence
$Q_1\in\cC_0(X)\subset\cJ_1$ and
\begin{equation}\label{eq:Krein3-9}
Q_1B_{12}\in\cJ_1.
\end{equation}
From (\ref{eq:Krein3-7})--(\ref{eq:Krein3-9}) we obtain that
$B_{12}=\cP A^{-1}\cQ\in\cJ_1$. By the definition of $\cL_1$, the operator
$A^{-1}$ belongs to $\cL_1$. Part (a) is proved.

(b) This statement follows from (a) with $\cQ=I-\cP$ in place of $\cP$.

(c) Since $\cP A\cQ\in\cJ_1\subset\cC_\infty(X)$,
$\cQ A\cP\in\cJ_2\subset\cC_\infty(X)$, and the operator $A$ is invertible,
we conclude that
\[
A_1=A-\cP A\cQ-\cQ A\cP=\cP A\cP+\cQ A\cQ
\]
is Fredholm. If $R$ is its regularizer, then $\cP R\cP|\im\cP$ is a regularizer
of $\cP A\cP|\im \cP$ on $\cP X$ and $\cQ R\cQ|\im\cQ$ is a regularizer of
$\cQ A\cQ|\im\cQ$ on $\cQ X$. Thus $\cP A\cP|\im \cP$ and $\cQ A\cQ|\im\cQ$
are Fredholm. Now statement (c) follows from parts (a) and (b).
\end{proof}
\begin{remark}
A minor modification of the proof of part (a) shows that if $A$ is invertible
on $X$, $\cJ_1=\{0\}$, and $\cP A\cP|\im\cP$ is invertible on $\cP X$, then
$A^{-1}\in\cL_1$.
\end{remark}
\begin{remark}
If $H$ is a separable Hilbert space and $\cJ$ is any two-sided ideal of
$\cB(H)$ such that $\cJ\ne\{0\}$ and $\cJ\ne\cB(H)$, then
$\cC_0(H)\subset\cJ\subset\cC_\infty(H)$
by Calkin's theorem (see, e.g. \cite[Chap. 3, Theorem~1.1]{GK69}).
Hence the statement of the above theorem can be simplified for separable
Hilbert spaces.
\end{remark}
\section{Auxiliary results}\label{sec:aux}
\subsection{Some facts on Besov spaces}
We start with the following well known facts.
\begin{lemma}\label{le:density-Besov}
If $1\le p<\infty$ and $0<\alpha\le 1$, then $\cP$ is dense in $B_p^\alpha$.
\end{lemma}
\begin{lemma}\label{le:embedding-Besov}
If $1< p<\infty$ and $1/p<\alpha\le 1$, then $B_p^\alpha\subset C$.
\end{lemma}
These lemmas are proved, for instance, in \cite[Sections~3.5.1 and 3.5.5]{ST87}.

The following fact is certainly known to specialists (see, e.g.
\cite[p. 735]{Pe03}). We give its proof for the convenience of the reader.
\begin{lemma}\label{le:A}
Let $1\le p<\infty$ and $0<\alpha<1$. If $a,b\in L^\infty\cap B_p^\alpha$, then
\[
\|ab\|_{B_p^\alpha}\le
\|a\|_{L^\infty}\|b\|_{B_p^\alpha}+
\|a\|_{B_p^\alpha}\|b\|_{L^\infty},
\]
that is, $L^\infty\cap B_p^\alpha$ is a Banach algebra under the
quasi-submultiplicative norm
\[
\|a\|_{L^\infty\cap B_p^\alpha}:=\|a\|_{L^\infty}+\|a\|_{B_p^\alpha}.
\]
\end{lemma}
\begin{proof}
It is easy to see that
\[
\omega_{L^p}^1(ab,s)\le
\|a\|_{L^\infty}\omega_{L^p}^1(b,s)+\omega_{L^p}^1(a,s)\|b\|_{L^\infty}
\quad (s\ge 0).
\]
From this inequality and
$\|ab\|_{L^p}\le\|a\|_{L^\infty}\|b\|_{L^p}+\|a\|_{L^p}\|b\|_{L^\infty}$
it follows that
\[
\begin{split}
\|ab\|_{B_p^\alpha}
&\le
\|a\|_{L^\infty}\|b\|_{L^p}+\|a\|_{L^p}\|b\|_{L^\infty}
\\
& +
\|a\|_{L^\infty}
\left(\int_0^{2\pi}\big[s^{-\alpha}\omega_{L^p}^1(b,s)\big]^p\frac{ds}{s}\right)^{1/p}
+
\|b\|_{L^\infty}
\left(\int_0^{2\pi}\big[s^{-\alpha}\omega_{L^p}^1(a,s)\big]^p\frac{ds}{s}\right)^{1/p}
\\
&=
\|a\|_{L^\infty}\|b\|_{B_p^\alpha}+
\|a\|_{B_p^\alpha}\|b\|_{L^\infty}.
\end{split}
\]
Now the fact that $\|\cdot\|_{L^\infty\cap B_p^\alpha}$ is a quasi-submultiplicative
norm is obvious.
\end{proof}
\subsection{The algebra of multiplication operators}
For $\delta,\mu\in\R$, we denote by $\ell_2^{\delta,\mu}$ the set of all
sequences $\varphi=\{\varphi_j\}_{j\in\Z}$ such that
\[
\|\varphi\|_{\ell_2^{\delta,\mu}}^2:=
\sum_{j=-\infty}^{-1}|\varphi_j|^2(|j|+1)^{2\delta}
+
\sum_{j=0}^\infty|\varphi_j|^2(j+1)^{2\mu}
<\infty.
\]
It is clear that $\ell_2^{\delta,\mu}$ is a Hilbert space.
Let $\ell^0$ denote the collection of all sequences from $\ell_2^{\delta,\mu}$
with finite support. For a function $a\in L^1$, define $M(a)$ on $\ell^0$ by
\[
M(a):\{\varphi_j\}_{j\in\Z}\mapsto\Big\{\sum_{k\in\Z}a_{j-k}\varphi_k\Big\}_j.
\]
If
\[
\sup\big\{
\|M(a)\varphi\|_{\ell_2^{\delta,\mu}}/\|\varphi\|_{\ell_2^{\delta,\mu}}\ :
\ \varphi\in\ell^0,\ \varphi\ne 0
\big\}<\infty,
\]
then $M(a)$ can be extended to a bounded operator on $\ell_2^{\delta,\mu}$.
In this case we call $M(a)$ the \textit{multiplication operator} with symbol $a$.
The following basic properties of multiplication operators on $\ell_2^{\delta,\mu}$
can be proved in the same way as in \cite[Sections~2.5 and~6.2]{BS06}
(see also \cite{V77}).
\begin{theorem}\label{th:M}
Let $\delta,\mu\in\R$ and $a\in L^1$.
\begin{enumerate}
\item[(a)]
If $M(a)\in\cB(\ell_2^{\delta,\mu})$, then the adjoint of $M(a)$ equals
$M(\overline{a})\in\cB(\ell_2^{-\delta,-\mu})$ and
\[
\|M(a)\|_{\cB(\ell_2^{\delta,\mu})}
=
\|M(\overline{a})\|_{\cB(\ell_2^{\delta,\mu})}
=
\|M(a)\|_{\cB(\ell_2^{-\delta,-\mu})}
=
\|M(\overline{a})\|_{\cB(\ell_2^{-\delta,-\mu})}.
\]
\item[(b)]
If $M(a)\in\cB(\ell_2^{\delta,\mu})$, then $a\in L^\infty$ and
\[
\|a\|_{L^\infty}
=
\|M(a)\|_{\cB(\ell_2^{0,0})}
\le
\|M(a)\|_{\cB(\ell_2^{\delta,\mu})}.
\]

\item[(c)]
The set $\{a\in L^1: M(a)\in\cB(\ell_2^{\delta,\mu})\}$ is a Banach algebra
under the Banach algebra norm $\|a\|:=\|M(a)\|_{\cB(\ell_2^{\delta,\mu})}$.
\end{enumerate}
\end{theorem}
\begin{lemma}\label{le:SR}
Let $0<\delta<\infty$. For every $k\in\Z$,
\[
\lim_{m\to+\infty}\|M(\chi_{km})\|_{\cB(\ell_2^{\delta,-\delta})}^{1/m}\le 1.
\]
\end{lemma}
\begin{proof}
It is easy to see that $M(\chi_{km})$ is the shift operator
$\{\varphi_j\}_{j\in\Z}\mapsto\{\varphi_{j-km}\}_{j\in\Z}$.
Let $\varphi\in\ell_2^{\delta,-\delta}$ and suppose $k$ is positive.
Without loss of generality assume $km\ge 1$. Then
\begin{equation}\label{eq:SR-1}
\begin{split}
\|M(\chi_{km})\varphi\|_{\ell_2^{\delta,-\delta}}^2
&=
\sum_{j=-\infty}^{-1}|\varphi_{j-km}|^2(|j|+1)^{2\delta}
+
\sum_{j=0}^\infty|\varphi_{j-km}|^2(|j|+1)^{-2\delta}
\\
&=
\sum_{i=-\infty}^{-1+km}|\varphi_i|^2(|i+km|+1)^{2\delta}
+
\sum_{i=km}^\infty|\varphi_i|^2(|i+km|+1)^{-2\delta}
\\
&\le
\sup_{i\le -1}\left(\frac{|i+km|+1}{|i|+1}\right)^{2\delta}
\sum_{i=-\infty}^{-1}|\varphi_i|^2(|i|+1)^{2\delta}
\\
&\quad+
\max_{0\le i\le km-1}\Big[(|i+km|+1)^{2\delta}(|i|+1)^{2\delta}\Big]
\sum_{i=0}^{km-1}|\varphi_i|^2(|i|+1)^{-2\delta}
\\
&\quad+
\sup_{i\ge km}\left(\frac{|i+km|+1}{|i|+1}\right)^{-2\delta}
\sum_{i=km}^\infty|\varphi_i|^2(|i|+1)^{-2\delta}
\\
&\le
\big(S_1(k,m)+S_2(k,m)+S_3(k,m)\big)
\|\varphi\|_{\ell_2^{\delta,-\delta}}^2
\end{split}
\end{equation}
where
\[
\begin{split}
S_1(k,m) &:=
\sup_{i\le -1}\left(\frac{|i+km|+1}{|i|+1}\right)^{2\delta},
\\
S_2(k,m) &:=
\max_{0\le i\le km-1}\Big[(|i+km|+1)^{2\delta}(|i|+1)^{2\delta}\Big],
\\
S_3(k,m) &:=
\sup_{i\ge km}\left(\frac{|i+km|+1}{|i|+1}\right)^{-2\delta}.
\end{split}
\]
If $i\le -1$, then
\[
\frac{|i+km|+1}{|i|+1}
\le
1+\frac{km}{|i|+1}
\le
1+\frac{km}{2}
\le
2km.
\]
Hence
\begin{equation}\label{eq:SR-2}
S_1(k,m)\le (2km)^{2\delta}.
\end{equation}
If $0\le i\le km-1$, then
\[
(|i+km|+1)(|i|+1)=(i+km+1)(i+1)\le 2(km)^2,
\]
whence
\begin{equation}\label{eq:SR-3}
S_2(k,m)\le 2^{2\delta}(km)^{4\delta}.
\end{equation}
If $i\ge km$, then
\[
\frac{|i+km|+1}{|i|+1}=\frac{i+km+1}{i+1}\ge 1.
\]
Therefore,
\begin{equation}\label{eq:SR-4}
S_3(k,m)\le 1.
\end{equation}
Combining (\ref{eq:SR-1})--(\ref{eq:SR-4}), we get for
$\varphi\in\ell_2^{\delta,-\delta}$,
\[
\|M(\chi_{km})\varphi\|_{\ell_2^{\delta,-\delta}}^2
\le
\big((2km)^{2\delta}+2^{2\delta}(km)^{4\delta}+1\big)
\|\varphi\|_{\ell_2^{\delta,-\delta}}^2.
\]
Thus,
\[
\lim_{m\to+\infty}\|M(\chi_{km})\|_{\cB(\ell_2^{\delta,-\delta})}^{1/m}
\le
\lim_{m\to+\infty}\big((2km)^{2\delta}+2^{2\delta}(km)^{4\delta}+1\big)^{1/(2m)}=1.
\]
If $k$ is negative, then the proof is analogous.
\end{proof}
Denote by $P$ the operator given on $\ell_2^{\delta,\mu}$ by
$(P\varphi)_j=\varphi_j$ if $j\ge 0$ and $(P\varphi)_j=0$ if $j<0$.
Let $Q:=I-P$. It is easy to see that $P^2=P$, $Q^2=Q$, and $\|P\|=\|Q\|=1$.
\begin{lemma}\label{le:simple}
Suppose $\mu\in\R$ and $a\in L^1$. Then
\[
\|PM(a)P\|_{\cB(\ell_2^{\mu,-\mu})}
\le
\|PM(a)P\|_{\cB(\ell_2^{-\mu,-\mu})},
\quad
\|QM(a)Q\|_{\cB(\ell_2^{\mu,-\mu})}
\le
\|QM(a)Q\|_{\cB(\ell_2^{\mu,\mu})}.
\]
\end{lemma}
\begin{proof}
This statement follows from the definition of $P$ and $Q$.
\end{proof}
The Besov space $B_\infty^\gamma$, $0<\gamma<1$, is nothing else than the
H\"older space $C^\gamma$, $0<\gamma<1$, defined as the set of all functions
$f\in C$ such that
\[
\|f\|_{C^\gamma}:=\|f\|_C+\sup_{0<s\le 2\pi}\frac{\omega_C^1(f,s)}{s^\gamma}<\infty.
\]
It is easy to see that $\|\cdot\|_{C^\gamma}$ is a Banach algebra norm.
\begin{theorem}[Verbitsky \cite{V77}]
\label{th:Verbitsky}
If $|\mu|<\gamma<1$, then there exists a positive constant
$L_{\gamma,\mu}$ depending only on $\gamma$ and $\mu$ such that
$\|M(a)\|_{\cB(\ell_2^{\mu,\mu})}\le L_{\gamma,\mu}\|a\|_{C^\gamma}$
for all $a\in C^\gamma$.
\end{theorem}
\subsection{Peller's theorems on Hankel-type operators}
For $\gamma\in(0,1)$, denote by $c^\gamma$ the closure of the set of all
Laurent polynomials in the norm of $C^\gamma$.

The following result is a corollary of Peller's theorems
\cite[Chap.~6, Theorem~8.1 and 8.2]{Pe03}.
It deals with the boundedness and compactness of $QM(a)P$ and $PM(a)Q$ on the
spaces $\ell_2^{\gamma/2,-\gamma/2}$ and $\ell_2^{-\gamma/2,\gamma/2}$, respectively.
\begin{theorem}\label{th:Peller1}
If $0<\gamma<1$ and $a\in L^\infty$, then
\begin{eqnarray}
\label{eq:Peller1-1}
&&
QM(a)P\in\cB(\ell_2^{\gamma/2,-\gamma/2})
\Longleftrightarrow
Qa\in C^\gamma
\Longleftrightarrow
a\in C^\gamma+H^\infty,
\\
\label{eq:Peller1-2}
&&
QM(a)P\in\cC_\infty(\ell_2^{\gamma/2,-\gamma/2})
\Longleftrightarrow
Qa\in c^\gamma
\Longleftrightarrow
a\in c^\gamma+H^\infty,
\\
\nonumber
&&
PM(a)Q\in\cB(\ell_2^{-\gamma/2,\gamma/2})
\Longleftrightarrow
Pa\in C^\gamma
\Longleftrightarrow
a\in C^\gamma+\overline{H^\infty},
\\
\nonumber
&&
PM(a)Q\in\cC_\infty(\ell_2^{-\gamma/2,\gamma/2})
\Longleftrightarrow
Pa\in c^\gamma
\Longleftrightarrow
a\in c^\gamma+\overline{H^\infty},
\end{eqnarray}
and there exist positive constants $c_1$ and $c_2$ depending only on $\gamma$
such that
\begin{eqnarray}
\label{eq:Peller1-3}
&&
c_1\|Qa\|_{C^\gamma}\le
\|QM(a)P\|_{\cB(\ell_2^{\gamma/2,-\gamma/2})}
\le
c_2\|Qa\|_{C^\gamma},
\\
\nonumber
&&
c_1\|Pa\|_{C^\gamma}\le
\|PM(a)Q\|_{\cB(\ell_2^{-\gamma/2,\gamma/2})}
\le
c_2\|Pa\|_{C^\gamma}.
\end{eqnarray}
\end{theorem}
The following theorem is a particular case of Peller's description of
generalized Hankel matrices belonging to the Schatten-von Neumann classes
(see \cite[Chap.~6, Theorem~8.9]{Pe03}).
\begin{theorem}\label{th:Peller2}
Let $1\le p,q<\infty$ and $\delta,\mu\in\R$. Suppose $a\in L^\infty$.
\begin{enumerate}
\item[(a)]
If  $0<1/p+\delta+\mu\le 1$ and
\begin{equation}\label{eq:Peller2-1}
\min\{\delta,\mu\}>\max\{-1/2,-1/p\},
\end{equation}
then
\[
QM(a)P\in\cC_p(\ell_2^{\delta,-\mu})
\Longleftrightarrow
Qa\in B_p^{1/p+\delta+\mu}
\Longleftrightarrow
a\in B_p^{1/p+\delta+\mu}+H^\infty,
\]
and there exist positive constants $c_1$ and $c_2$ depending only on
$\delta,\mu$, and $p$ such that
\[
c_1\|Qa\|_{B_p^{1/p+\delta+\mu}}
\le
\|QM(a)P\|_{\cC_p(\ell_2^{\delta,-\mu})}
\le
c_2\|Qa\|_{B_p^{1/p+\delta+\mu}}.
\]

\item[(b)]
If  $0<1/q-\delta-\mu\le 1$ and
\begin{equation}\label{eq:Peller2-2}
\max\{\delta,\mu\}<\min\{1/2,1/q\},
\end{equation}
then
\[
PM(a)Q\in\cC_q(\ell_2^{\delta,-\mu})
\Longleftrightarrow
Pa\in B_q^{1/q-\delta-\mu}
\Longleftrightarrow
a\in B_q^{1/q-\delta-\mu}+\overline{H^\infty},
\]
and there exist positive constants $c_1$ and $c_2$ depending only on $\delta,\mu$,
and $q$ such that
\[
c_1\|Pa\|_{B_q^{1/q-\delta-\mu}}
\le
\|PM(a)Q\|_{\cC_q(\ell_2^{\delta,-\mu})}
\le
c_2\|Pa\|_{B_q^{1/q-\delta-\mu}}.
\]
\end{enumerate}
\end{theorem}
\begin{remark}
The restrictions $1\le p,q<\infty$ and $0<1/p+\delta+\mu\le 1$
(respectively, $0<1/q-\delta-\mu\le 1$) are not essential. The theorem is also
true for $p,q\in(0,1)$ and large $\delta$ and $\mu$ (resp. large $-\delta$ and
$-\mu$). We imposed these restrictions just to keep the presentation in the
setting of Besov spaces $B_p^\alpha$ with $1\le p<\infty$ and $0<\alpha\le 1$,
which is sufficient for our purposes.

On the other hand, hypotheses (\ref{eq:Peller2-1}) and (\ref{eq:Peller2-2}) are
essential because without these hypotheses the theorem is not true (see the remark
on p. 291 of \cite{Pe03}).
\end{remark}
\subsection{Toeplitz operators with antianalytic symbols on analytic Besov spaces}
The following result is due to Peller and Khrushchev \cite{PK82}.
Its proof is contained in \cite[Proposition 10.23]{BS06}.
\begin{lemma}\label{le:L}
Let $1<p<\infty$ and $0<\alpha<1$. If $g\in H^\infty$, then the Toeplitz
operator $T(\overline{g}):\varphi\mapsto P(\overline{g}\varphi)$ is bounded
on $PB_p^\alpha$ and there exists a positive constant $L_{p,\alpha}$ depending
only on $p$ and $\alpha$ such that
\[
\|T(\overline{g})\|_{\cB(PB_p^\alpha)}\le L_{p,\alpha}\|g\|_{L^\infty}.
\]
\end{lemma}
\section{The Banach algebras $C^\gamma+H^\infty$ and $c^\gamma+H^\infty$}
\label{sec:CH-and-cH}
\subsection{The sets $C^\gamma+H^\infty$ and $c^\gamma+H^\infty$ are Banach algebras}
The results of this section may look curious at the first glance. However, they
are important pieces of the proof of Theorem~\ref{th:invertibility}. The material
of this section is taken from \cite[Chap.~4]{BS83}.
\begin{lemma}\label{le:identification}
If $0<\gamma<1$ and $a\in L^1$, then
\[
M(a)\in\cB(\ell_2^{\gamma/2,-\gamma/2})
\Longleftrightarrow
a\in C^\gamma+H^\infty.
\]
\end{lemma}
\begin{proof}
If $M(a)\in\cB(\ell_2^{\gamma/2,-\gamma/2})$, then $a\in L^\infty$ in view
of Theorem~\ref{th:M}(b). It is clear that $QM(a)P\in\cB(\ell_2^{\gamma/2,-\gamma/2})$.
By (\ref{eq:Peller1-1}), $a\in C^\gamma+H^\infty$. The necessity portion is
proved.

Let us prove the sufficiency part. If $a\in C^\gamma+H^\infty$, then there exist
$c\in C^\gamma$ and $h\in H^\infty$ such that $a=c+h$. We represent $M(a)$
as
\[
M(a)=PM(c)P+PM(h)P+QM(a)P+PM(a)Q+QM(c)Q+QM(h)Q
\]
and show that all terms on the right-hand side are bounded on
$\ell_2^{\gamma/2,-\gamma/2}$.

By Theorem~\ref{th:Verbitsky}, $PM(c)P\in\cB(\ell_2^{-\gamma/2,-\gamma/2})$
and $QM(c)Q\in\cB(\ell_2^{\gamma/2,\gamma/2})$. Thus, in view of
Lemma~\ref{le:simple}, $PM(c)P$ and $QM(c)Q$ are bounded on
$\ell_2^{\gamma/2,-\gamma/2}$.

Lemma~\ref{le:L} yields that the operators $T(\overline{h})$ and $T(\widetilde{h})$
are bounded on $PB_2^{\gamma/2}$. Hence the operators $PM(\overline{h})P$ and
$QM(h)Q$ are bounded on the Besov space $B_2^{\gamma/2}$. It is well known
that $f\in B_2^{\gamma/2}$ if and only if its sequence of the Fourier coefficients
belongs to $\ell_2^{\gamma/2,\gamma/2}$ and that the corresponding norms are equivalent.
Thus $PM(\overline{h})P$ and $QM(h)Q$ are bounded on $\ell_2^{\gamma/2,\gamma/2}$.
By Theorem~\ref{th:M}(a), the operator $PM(h)P$ is the adjoint of $PM(\overline{h})P$
and $PM(h)P\in\cB(\ell_2^{-\gamma/2,-\gamma/2})$. From Lemma~\ref{le:simple} it
follows that the operators $PM(h)P$ and $QM(h)Q$ are bounded on
$\ell_2^{\gamma/2,-\gamma/2}$.

Let $\varphi\in\ell_2^{\gamma/2,-\gamma/2}$. Then, taking into account that
$\|M(a)\|_{\cB(\ell_2^{0,0})}=\|a\|_{L^\infty}$, we have
\[
\begin{split}
\|PM(a)Q\varphi\|_{\ell_2^{\gamma/2,-\gamma/2}}
&=
\|PM(a)Q\varphi\|_{\ell_2^{0,-\gamma/2}}
\le
\|PM(a)Q\varphi\|_{\ell_2^{0,0}}
\\
&\le
\|a\|_{L^\infty}\|Q\varphi\|_{\ell_2^{0,0}}
=
\|a\|_{L^\infty}\|Q\varphi\|_{\ell_2^{0,-\gamma/2}}
\le
\|a\|_{L^\infty}\|\varphi\|_{\ell_2^{\gamma/2,-\gamma/2}}.
\end{split}
\]
Thus $PM(a)Q$ is bounded on $\ell_2^{\gamma/2,-\gamma/2}$.
Finally, the operator $QM(a)P$ is bounded on
$\ell_2^{\gamma/2,-\gamma/2}$ in view of (\ref{eq:Peller1-1}).
\end{proof}
\begin{theorem}\label{th:CH-algebra}
If $0<\gamma<1$, then the sets $C^\gamma+H^\infty$ and $c^\gamma+H^\infty$
are Banach algebras under the Banach algebra norm
\[
\|a\|:=\|M(a)\|_{\cB(\ell_2^{\gamma/2,-\gamma/2})}.
\]
\end{theorem}
\begin{proof}
From the sufficiency portion of Lemma~\ref{le:identification} it follows that
\[
C^\gamma+H^\infty=\{a\in C^\gamma+H^\infty:M(a)\in\cB(\ell_2^{\gamma/2,-\gamma/2})\}.
\]
It is obvious that $C^\gamma+H^\infty$ is an algebra. Let $a_n$ be
a Cauchy sequence in this algebra. By Theorem~\ref{th:M}(c), there
exists an $a\in L^1$ such that
$\|M(a_n)-M(a)\|_{\cB(\ell_2^{\gamma/2,-\gamma/2})}\to 0$ as
$n\to\infty$ and $M(a)\in\cB(\ell_2^{\gamma/2,-\gamma/2})$. In
view of the necessity part of Lemma~\ref{le:identification}, $a\in
C^\gamma+H^\infty$. Hence $C^\gamma+H^\infty$ is closed. It is
clear that $\|\cdot\|$ is a Banach algebra norm.

Let $\cA:=\cB(\ell_2^{\gamma/2,-\gamma/2})$ and $\cL:=\{M(a)\in\cA:a\in C^\gamma+H^\infty\}$.
We have proved that $\cL$ is a closed subalgebra of $\cA$. From
Theorem~\ref{th:Krein3}(b) with $\cA$, $\cL$, and
$\cJ_2:=\cC_\infty(\ell_2^{\gamma/2,-\gamma/2})$ we obtain that
$\cL_2:=\{M(a)\in\cL:QM(a)P\in\cJ_2\}$ is a Banach algebra under the Banach algebra
norm
\[
\|M(a)\|_1:=
\|M(a)\|_{\cB(\ell_2^{\gamma/2,-\gamma/2})}+
\|QM(a)P\|_{\cB(\ell_2^{\gamma/2,-\gamma/2})}.
\]
It is clear that this norm is equivalent to $\|M(a)\|_{\cB(\ell_2^{\gamma/2,-\gamma/2})}$.
From (\ref{eq:Peller1-2}) we get
\[
\cL_2=\{M(a)\in\cB(\ell_2^{\gamma/2,-\gamma/2})\ :\ a\in c^\gamma+H^\infty\}.
\]
Hence $c^\gamma+H^\infty$ is a Banach algebra under the norm $\|\cdot\|$.
\end{proof}
\subsection{Structure of the Banach algebra $c^\gamma+H^\infty$}
In this subsection we clarify the structure of the algebra $c^\gamma+H^\infty$.
\begin{lemma}\label{le:H-in-cH}
The set $H^\infty$ is a closed subalgebra of $c^\gamma+H^\infty$.
\end{lemma}
\begin{proof}
It is obvious that $H^\infty$ is a subalgebra of $c^\gamma+H^\infty$. Let
$h_n\in H^\infty$ be a Cauchy sequence in the norm of $c^\gamma+H^\infty$.
Since $c^\gamma+H^\infty$ is closed, there is a function $h\in c^\gamma+H^\infty$
such that $\|M(h_n)-M(h)\|_{\cB(\ell_2^{\gamma/2,-\gamma/2})}\to 0$ as
$n\to\infty$. From this fact and Theorem~\ref{th:M}(b) we deduce that
$\|h_n-h\|_{L^\infty}\to 0$ as $n\to\infty$. Because $H^\infty$ is closed
in $L^\infty$, we conclude that $h\in H^\infty$, that is, $H^\infty$
is closed in $c^\gamma+H^\infty$.
\end{proof}
\begin{lemma}\label{le:structure-cH}
We have $c^\gamma+H^\infty=\mathrm{clos}_{c^\gamma+H^\infty}\mathrm{span}
\{\psi\overline{\chi_n}:\psi\in H^\infty,n\in\Z_+\}$.
\end{lemma}
\begin{proof}
Put
$A:=\mathrm{clos}_{c^\gamma+H^\infty}\mathrm{span}
\{\psi\overline{\chi_n}:\psi\in H^\infty,n\in\Z_+\}$.
It is obvious that $A\subset c^\gamma+H^\infty$. Let us show the reverse
inclusion. If $f\in c^\gamma+H^\infty$, then there exist $c\in c^\gamma$ and
$h\in H^\infty$ such that $f=c+h$. By Privalov's theorem, the projections $P$
and $Q$ are bounded on $C^\gamma$, $0<\gamma<1$. Hence $Pc\in H^\infty$. It is
clear that $Pc+h\in H^\infty\subset A$. By the definition of $c^\gamma$, there
exists a sequence of Laurent polynomials $p_m$ such that
\begin{equation}\label{eq:structure-cH-1}
\|c-p_m\|_{C^\gamma}\to 0\quad\mbox{as}\quad m\to\infty.
\end{equation}
Obviously, $PM(Qc-Qp_m)Q=0$. Hence
\begin{equation}\label{eq:structure-cH-2}
M(Qc)-M(Qp_m)
=
PM(Qc-Qp_m)P
+
QM(Qc-Qp_m)P
+
QM(Qc-Qp_m)Q.
\end{equation}
From Theorem~\ref{th:Verbitsky} and Lemma~\ref{le:simple} it follows that
\begin{equation}\label{eq:structure-cH-3}
\begin{split}
\|PM(Qc-Qp_m)P\|_{\cB(\ell_2^{\gamma/2,-\gamma/2})}
&\le
\|PM(Qc-Qp_m)P\|_{\cB(\ell_2^{-\gamma/2,-\gamma/2})}
\\
&\le
\|M(Qc-Qp_m)\|_{\cB(\ell_2^{-\gamma/2,-\gamma/2})}
\\
&\le
L_{\gamma,-\gamma/2}\|Q\|_{\cB(C^\gamma)}\|c-p_m\|_{C^\gamma}
\end{split}
\end{equation}
and similarly,
\begin{equation}\label{eq:structure-cH-4}
\|QM(Qc-Qp_m)Q\|_{\cB(\ell_2^{\gamma/2,-\gamma/2})}
\le
L_{\gamma,\gamma/2}\|Q\|_{\cB(C^\gamma)}\|c-p_m\|_{C^\gamma}.
\end{equation}
On the other hand, from (\ref{eq:Peller1-3}) we see that
\begin{equation}\label{eq:structure-cH-5}
\|QM(Qc-Qp_m)P\|_{\cB(\ell_2^{\gamma/2,-\gamma/2})}
\le
c_2\|Q(Qc-Qp_m)\|_{C^\gamma}
\le
c_2\|Q\|_{\cB(C^\gamma)}\|c-p_m\|_{C^\gamma}.
\end{equation}
Combining (\ref{eq:structure-cH-1})--(\ref{eq:structure-cH-5}), we get
\[
\|M(Qc)-M(Qp_m)\|_{\cB(\ell_2^{\gamma/2,-\gamma/2})}
\le
\mathrm{const}
\|c-p_m\|_{C^\gamma}\to 0
\quad\mbox{as}\quad m\to\infty.
\]
It is clear that
$Qp_m\in\mathrm{span}\{\psi\overline{\chi_n}:\psi\in H^\infty,n\in\Z_+\}$
and therefore $Qc\in A$. Thus $f\in A$ and $c^\gamma+H^\infty\subset A$.
\end{proof}
\subsection{Invertibility in the Banach algebra $c^\gamma+H^\infty$}
In this subsection we will show that $c^\gamma+H^\infty$ is inverse closed in
$C+H^\infty$. From \cite[Proposition~6.36, Corollary~6.38]{D98} we get the
following description of the maximal ideal space of $C+H^\infty$.
\begin{lemma}\label{le:maximal-CH}
We have $M_{C+H^\infty}=
\{m\in M_{H^\infty}:|\widehat{\chi_k}(m)|=1\mbox{ for }k\in\Z_+\}$.
\end{lemma}
The next statement shows that the maximal ideal spaces of $c^\gamma+H^\infty$
and $C+H^\infty$ coincide.
\begin{lemma}\label{le:maximal-cH}
We have
$M_{c^\gamma+H^\infty}=
\{m\in M_{H^\infty}:|\widehat{\chi_k}(m)|=1\mbox{ for }k\in\Z_+\}$.
\end{lemma}
\begin{proof}
The proof is developed by analogy with \cite[Proposition~6.37]{D98}.
By Lemma~\ref{le:H-in-cH}, $H^\infty$ is a closed subalgebra of
$c^\gamma+H^\infty$. If $f$ is a multiplicative linear functional on
$c^\gamma+H^\infty$, then $f|H^\infty$ is a multiplicative linear
functional on $H^\infty$. Hence $\Psi:f\mapsto f|H^\infty$
defines a continuous mapping from $M_{c^\gamma+H^\infty}$ into
$M_{H^\infty}$. If $f_1$ and $f_2$ are elements in $M_{c^\gamma+H^\infty}$
such that $\Psi(f_1)=\Psi(f_2)$, then for $k\in\Z_+$,
\[
f_1(\overline{\chi_k})=f_1(\chi_k^{-1})=f_1(\chi_k)^{-1}
=f_2(\chi_k)^{-1}=f_2(\chi_k^{-1})=f_2(\overline{\chi_k}).
\]
From this equality and Lemma~\ref{le:structure-cH} we conclude that $f_1=f_2$.
Therefore, $\Psi$ is a homeomorphism of $M_{c^\gamma+H^\infty}$ into
$M_{H^\infty}$. Moreover, for every $f\in M_{c^\gamma+H^\infty}$ and
$k\in\Z_+$, by multiplicativity of $f$ and Lemma~\ref{le:SR},
\[
|f(\chi_k)|
\le
\lim_{m\to+\infty}\|\chi_k^m\|_{c^\gamma+H^\infty}^{1/m}
=
\lim_{m\to+\infty}\|M(\chi_{km})\|_{\cB(\ell_2^{\gamma/2,-\gamma/2})}^{1/m}
\le 1
\]
and similarly
\[
1/|f(\chi_k)|=|f(\chi_{-k})|\le 1.
\]
That is, $|f(\chi_k)|=1$ for every $k\in\Z_+$. Therefore, the range of $\Psi$
is contained in
\[
\big\{m\in M_{H^\infty}\ :\ |\widehat{\chi_k}(m)|=1\mbox{ for }k\in\Z_+\big\},
\]
and only the reverse containment remains.

Let $m$ be a point in $M_{H^\infty}$ such that $|\widehat{\chi_k}(m)|=1$ for
$k\in\Z_+$. If we define $g$ on
$\mathrm{span}\{\psi\overline{\chi_n}:\psi\in H^\infty, n\in\Z_+\}$
such that $g(\psi\overline{\chi_k})=\widehat{\psi}(m)\widehat{\chi_k}(m)$,
then $g$ can be easily shown to be multiplicative. The inequality
\[
\begin{split}
|g(\psi\overline{\chi_k})|
&=
|\widehat{\psi}(m)|\,|\widehat{\chi_k}(m)|
=
|\widehat{\psi}(m)|
\le
\|\psi\|_{H^\infty}
\\
&=
\|\psi\overline{\chi_k}\|_{L^\infty}
\le
\|M(\psi\overline{\chi_k})\|_{\cB(\ell_2^{\gamma/2,-\gamma/2})}
=
\|\psi\overline{\chi_k}\|_{c^\gamma+H^\infty}
\end{split}
\]
(recall Theorems~\ref{th:M}(b) and \ref{th:CH-algebra}) shows that $g$ can be
extended to a multiplicative linear functional on
$\mathrm{clos}_{c^\gamma+H^\infty}\mathrm{span}
\{\psi\overline{\chi_n}:\psi\in H^\infty, n\in\Z_+\}$.
But the latter set coincides with $c^\gamma+H^\infty$ by
Lemma~\ref{le:structure-cH}. Obviously, $\Psi(g)=m$ and thus
$M_{c^\gamma+H^\infty}$ is homeomorphic to
$\{m\in H^\infty:|\widehat{\chi_k}(m)|=1\mbox{ for }k\in\Z_+\}$.
\end{proof}
Combining Lemmas~\ref{le:maximal-CH} and \ref{le:maximal-cH} with Gelfand's
theorem, we get the following.
\begin{theorem}\label{th:invertibility-cH}
If $0<\gamma<1$ and $a\in c^\gamma+H^\infty$, then
\[
a\in G(c^\gamma+H^\infty)\Longleftrightarrow a\in G(C+H^\infty).
\]
\end{theorem}
\section{Generalized Krein algebras}
\label{sec:gen-Krein-algebras}
\subsection{The sets $K_{p,0}^{\alpha,0}$, $K_{0,q}^{0,\beta}$, and
$K_{p,q}^{\alpha,\beta}$ are Banach algebras}
Since the projections $P$ and $Q$ are bounded on the Besov spaces $B_p^\alpha$ and
$B_q^\beta$, it is clear that $K_{p,0}^{\alpha,0}$, $K_{0,q}^{0,\beta}$,
and $K_{p,q}^{\alpha,\beta}$ are Banach spaces under the norms
(\ref{eq:def-gen-Krein1}),
(\ref{eq:def-gen-Krein2}), and
(\ref{eq:def-gen-Krein3}), respectively.
In this subsection we show that these norms are quasi-submultiplicative.
\begin{lemma}\label{le:PQ}
Let $1<p,q<\infty$ and $0<\alpha,\beta<1$.
\begin{enumerate}
\item[(a)]
If $\alpha>1/p$, then the projections $P$ and $Q$ are bounded on $K_{p,0}^{\alpha,0}$.
\item[(b)]
If $\beta>1/q$, then the projections $P$ and $Q$ are bounded on $K_{0,q}^{0,\beta}$.
\item[(c)]
If $\alpha>1/p$ or $\beta>1/q$, then the projections $P$ and $Q$ are
bounded on $K_{p,q}^{\alpha,\beta}$.
\end{enumerate}
\end{lemma}
\begin{proof}
(a) By Lemma~\ref{le:embedding-Besov}, there exists a constant $C>0$ such that
$\|f\|_{L^\infty}\le C\|f\|_{B_p^\alpha}$ for all $f\in B_p^\alpha$. If
$a\in K_{p,0}^{\alpha,0}$, then
\[
\|Qa\|_{K_{p,0}^{\alpha,0}}
=
\|Qa\|_{L^\infty}+\|Qa\|_{B_p^\alpha}
\le
(C+1)\|Qa\|_{B_p^\alpha}
\le
(C+1)\|a\|_{K_{p,0}^{\alpha,0}}.
\]
The boundedness of $P=I-Q$ is now obvious.

(b) The proof is analogous.

(c) The arguments of the proof of (a) or (b) apply.
\end{proof}
Now we are in a position to prove our first main result.
\begin{proof}[Proof of Theorem~{\rm\ref{th:algebra}}]
Let us prove only part (c). The proofs of parts (a) and (b) are similar
(and even a little simpler).

For $\alpha=1/p$ and $\beta=1/q$, the statement is contained in
Corollary~\ref{co:Krein-representation}. Hence we can assume that $\alpha> 1/p$ or
$\beta>1/q$.
Then the projections $P$ and $Q$ are bounded on $K_{p,q}^{\alpha,\beta}$
by Lemma~\ref{le:PQ}(c).
For brevity we will omit the subscript in $\|\cdot\|_{K_{p,q}^{\alpha,\beta}}$.
If $a,b\in K_{p,q}^{\alpha,\beta}$, then we have
\begin{equation}\label{eq:algebra-1}
\begin{split}
\|ab\|
&\le
\|PaPb\|+\|QaQb\|+\|PaQb\|+\|QaPb\|
\\
&=
\|PaPb\|_{L^\infty}+\|P(PaPb)\|_{B_q^\beta}
+
\|QaQb\|_{L^\infty}+\|Q(QaQb)\|_{B_p^\alpha}
\\
&\quad
+\|PaQb\|_{L^\infty}+\|P(PaQb)\|_{B_q^\beta}+\|Q(PaQb)\|_{B_p^\alpha}
\\
&\quad
+\|QaPb\|_{L^\infty}+\|P(QaPb)\|_{B_q^\beta}+\|Q(QaPb)\|_{B_p^\alpha}.
\end{split}
\end{equation}
It is clear that
\begin{equation}\label{eq:algebra-2}
\begin{split}
&
\|PaQb\|_{L^\infty}+
\|QaQb\|_{L^\infty}+
\|PaQb\|_{L^\infty}+
\|QaPb\|_{L^\infty}
\\
&
\le
\|Pa\|\,\|Pb\|+\|Qa\|\,\|Qb\|+
\|Pa\|\,\|Qb\|+\|Qa\|\,\|Pb\|
\\
&
\le
(\|P\|+\|Q\|)^2\|a\|\,\|b\|.
\end{split}
\end{equation}
In view of Lemma~\ref{le:A},
\begin{equation}\label{eq:algebra-4}
\begin{split}
\|P(PaPb)\|_{B_q^\beta}
&=
\|PaPb\|_{B_q^\beta}
\\
&\le
\|Pa\|_{L^\infty}
\|Pb\|_{B_q^\beta}
+
\|Pa\|_{B_q^\beta}
\|Pb\|_{L^\infty}
\\
&\le
\|Pa\|\,\|b\|+\|a\|\,\|Pb\|
\\
&\le
2\|P\|\,\|a\|\,\|b\|
\end{split}
\end{equation}
and analogously
\begin{equation}\label{eq:algebra-5}
\|Q(QaQb)\|_{B_p^\alpha}
\le
2\|Q\|\,\|a\|\,\|b\|.
\end{equation}
By Lemma~\ref{le:L},
\begin{equation}\label{eq:algebra-6}
\begin{split}
\|P(PaQb)\|_{B_q^\beta}
&=
\|PM(Qb)Pa\|_{B_q^\beta}
\\
&=
\|T(Qb)Pa\|_{PB_q^\beta}
\\
&\le
\|T(Qb)\|_{\cB(PB_q^\beta)}\|Pa\|_{PB_q^\beta}
\\
&\le
L_{q,\beta}\|Qb\|_{L^\infty}\|Pa\|_{B_q^\beta}
\\
&\le
L_{q,\beta}\|Qb\|\,\|a\|
\\
&\le
L_{q,\beta}\|Q\|\,\|a\|\,\|b\|
\end{split}
\end{equation}
and analogously
\begin{equation}\label{eq:algebra-7}
\|P(QaPb)\|_{B_q^\beta}
\le
L_{q,\beta}\|Q\|\,\|a\|\,\|b\|.
\end{equation}
It is clear that the operator $J$ is an isometry on $B_p^\alpha$. Then applying
Lemma~\ref{le:L} again, we get
\begin{equation}\label{eq:algebra-8}
\begin{split}
\|Q(PaQb)\|_{B_p^\alpha}
&=
\|QM(Pa)Qb\|_{B_p^\alpha}
\\
&=
\|JT(\widetilde{Pa})JQb\|_{QB_p^\alpha}
\\
&\le
\|T(\widetilde{Pa})\|_{\cB(PB_p^\alpha)}\|JQb\|_{PB_p^\alpha}
\\
&\le
L_{p,\alpha}\|Pa\|_{L^\infty}\|Qb\|_{B_p^\alpha}
\\
&\le
L_{p,\alpha}\|Pa\|\,\|b\|
\\
&\le
L_{p,\alpha}\|P\|\,\|a\|\,\|b\|
\end{split}
\end{equation}
and similarly
\begin{equation}\label{eq:algebra-9}
\|Q(QaPb)\|_{B_p^\alpha}\le L_{p,\alpha}\|P\|\,\|a\|\,\|b\|.
\end{equation}
Combining (\ref{eq:algebra-1})--(\ref{eq:algebra-9}), we get
\[
\|ab\|\le C(p,q,\alpha,\beta)\|a\|\,\|b\|
\]
with $C(p,q,\alpha,\beta):=
(\|P\|+\|Q\|)^2 +2\|P\|+2\|Q\|+2L_{q,\beta}\|Q\|+2L_{p,\alpha}\|P\|$.
\end{proof}
\subsection{Invertibility in the Banach algebras
$K_{p,0}^{\alpha,0}$, $K_{0,q}^{0,\beta}$, and $K_{p,q}^{\alpha,\beta}$}
In this subsection we show that generalized Krein algebras are inverse closed
either in $C+H^\infty$ or in $C+\overline{H^\infty}$.

From Lemma~\ref{le:embedding-Besov} and (\ref{eq:embedding}) we immediately get
the following.
\begin{lemma}\label{le:embedding-Krein}
Let $1<p,q<\infty$ and $0<\alpha,\beta<1$.
\begin{enumerate}
\item[(a)]
If $\alpha\ge 1/p$, then
$K_{p,0}^{\alpha,0}\subset C+H^\infty$ and
$K_{p,q}^{\alpha,\beta}\subset C+H^\infty$.
\item[(b)]
If $\beta\ge 1/q$, then
$K_{0,q}^{0,\beta}\subset C+\overline{H^\infty}$ and
$K_{p,q}^{\alpha,\beta}\subset C+\overline{H^\infty}$.
\end{enumerate}
\end{lemma}
Now we prepare the proof of Theorem~\ref{th:invertibility}. We know from
Theorem~\ref{th:invertibility-cH} that $c^\gamma+H^\infty$ is inverse closed
in $C+H^\infty$. We show that the intersection of a generalized Krein
algebra $K_{p,q}^{\alpha,\beta}$ (or $K_{p,0}^{\alpha,0}$) with $c^\gamma+H^\infty$
is inverse closed in $c^\gamma+H^\infty$ (and thus in $C+H^\infty$) if $\alpha>1/p$.
\begin{lemma}\label{le:IA}
Let $0<\lambda\le 1$, $1<p,q<\infty$, $1/p+1/q=\lambda$, and $0<\gamma<\lambda-1/p$.
Suppose $K$ is either $K_{p,q}^{1/p+\gamma,1/q-\gamma}$ or $K_{p,0}^{1/p+\gamma,0}$.
\begin{enumerate}
\item[(a)]
The set $K\cap (c^\gamma+H^\infty)$ is a
Banach algebra under the quasi-sub\-multiplicative norm
\[
\|a\|:=
\|M(a)\|_{\cB(\ell_2^{\gamma/2,-\gamma/2})}+\|a\|_K.
\]

\item[(b)]
If $a\in K\cap (c^\gamma+H^\infty)$, then
\[
a\in G\big(K\cap(c^\gamma+H^\infty)\big)
\Longleftrightarrow
a\in G(c^\gamma+H^\infty).
\]
\end{enumerate}
\end{lemma}
\begin{proof}
(a) This statement follows from Theorem~\ref{th:CH-algebra} and
Theorem~\ref{th:algebra}(a), (c).

(b) Let $\cA:=\cB(\ell_2^{\gamma/2,-\gamma/2})$ and
$\cL:=\{M(a)\in\cA:a\in c^\gamma+H^\infty\}$. From Theorem~\ref{th:CH-algebra}
it follows that the mapping
\[
M:c^\gamma+H^\infty\to\cL,\quad a\mapsto M(a)
\]
is an isometry of Banach algebras. Hence
\begin{equation}\label{eq:IA-1}
a\in G(c^\gamma+H^\infty)\Longleftrightarrow M(a)\in G\cL.
\end{equation}
By Lemma~\ref{le:Schatten-properties}(b), (c) and
Theorems~\ref{th:Krein1} and~\ref{th:Krein2},
\[
\begin{split}
\cL_2 &:=\big\{M(a)\in\cL\ :\
QM(a)P\in\cC_p(\ell_2^{\gamma/2,-\gamma/2})
\big\},
\\
\cL_* &:=\big\{M(a)\in\cL\ :\
PM(a)Q\in\cC_q(\ell_2^{\gamma/2,-\gamma/2}),\
QM(a)P\in\cC_p(\ell_2^{\gamma/2,-\gamma/2})
\big\}
\end{split}
\]
are Banach algebras under the Banach algebra norms
\[
\begin{split}
\|M(a)\|_2 &:=
\|M(a)\|_{\cB(\ell_2^{\gamma/2,-\gamma/2})}+
\|QM(a)P\|_{\cC_p(\ell_2^{\gamma/2,-\gamma/2})},
\\
\|M(a)\|_* &:=
\|M(a)\|_{\cB(\ell_2^{\gamma/2,-\gamma/2})}+
\|PM(a)Q\|_{\cC_q(\ell_2^{\gamma/2,-\gamma/2})}+
\|QM(a)P\|_{\cC_p(\ell_2^{\gamma/2,-\gamma/2})},
\end{split}
\]
respectively.

If $M(a)\in G\cL$, then from (\ref{eq:IA-1}) and (\ref{eq:Peller1-2})
it follows that $QM(a)P$ and $QM(a^{-1})P$ belong to
$\cC_\infty(\ell_2^{\gamma/2,-\gamma/2})$. Thus,
\[
\begin{split}
Q-QM(a)Q\cdot QM(a^{-1})Q &=QM(a)P\cdot PM(a^{-1})Q\in
\cC_\infty(\ell_2^{\gamma/2,-\gamma/2}),
\\
Q-QM(a^{-1})Q\cdot QM(a)Q &=QM(a^{-1})P\cdot PM(a)Q\in
\cC_\infty(\ell_2^{\gamma/2,-\gamma/2}),
\end{split}
\]
that is, $QM(a^{-1})Q|\im Q$ is a regularizer of $QM(a)Q|\im Q$ on
$Q\ell_2^{\gamma/2,-\gamma/2}$. By Theorem~\ref{th:Krein3}(b), $M(a)\in G\cL_2$.
In view of Theorem~\ref{th:Krein3}(c), we conclude also that $M(a)\in G\cL_*$.
Thus,
\begin{equation}\label{eq:IA-2}
M(a)\in G\cL\Longleftrightarrow M(a)\in G\cL_2,
\quad
M(a)\in G\cL\Longleftrightarrow M(a)\in G\cL_*.
\end{equation}
Take $\mu=\delta=\gamma/2$. Then hypotheses (\ref{eq:Peller2-1}) and
(\ref{eq:Peller2-2}) are simultaneously satisfied. Due to Theorem~\ref{th:Peller2},
\[
\begin{split}
\cL_2 &=\big\{
M(a)\in\cA\ : a\in K_{p,0}^{1/p+\gamma,0}\cap (c^\gamma+H^\infty)
\big\},
\\
\cL_* &=\big\{
M(a)\in\cA\ : a\in K_{p,q}^{1/p+\gamma,1/q-\gamma}\cap (c^\gamma+H^\infty)
\big\},
\end{split}
\]
and the norms $\|PM(a)Q\|_{\cC_q(\ell_2^{\gamma/2,-\gamma/2})}$
and $\|QM(a)P\|_{\cC_p(\ell_2^{\gamma/2,-\gamma/2})}$ are equivalent to the norms
$\|Pa\|_{B_q^{1/q-\gamma}}$ and $\|Qa\|_{B_p^{1/p+\gamma}}$, respectively.
From Theorem~\ref{th:M}(b) we conclude that
$\|M(a)\|_{\cB(\ell_2^{\gamma/2,-\gamma/2})}$ is equivalent to
$\|a\|_{L^\infty}+\|M(a)\|_{\cB(\ell_2^{\gamma/2,-\gamma/2})}$. Thus,
the norms $\|M(a)\|_2$ and $\|M(a)\|_*$ are equivalent to the norms
\[
\begin{split}
\|M(a)\|_2' &:=
\|M(a)\|_{\cB(\ell_2^{\gamma/2,-\gamma/2})}+
\|a\|_{K_{p,0}^{1/p+\gamma,0}},
\\
\|M(a)\|_*' &:=
\|M(a)\|_{\cB(\ell_2^{\gamma/2,-\gamma/2})}+
\|a\|_{K_{p,q}^{1/p+\gamma,1/q-\gamma}},
\end{split}
\]
respectively.
It is clear that the mappings
\[
\begin{split}
&
M:K_{p,0}^{1/p+\gamma,0}\cap(c^\gamma+H^\infty)\to(\cL_2,\|\cdot\|_2'),
\quad
a\mapsto M(a),
\\
&
M:K_{p,q}^{1/p+\gamma,1/q-\gamma}\cap(c^\gamma+H^\infty)\to(\cL_*,\|\cdot\|_*'),
\quad
a\mapsto M(a)
\end{split}
\]
are isometries. Thus,
\begin{equation}\label{eq:IA-3}
\begin{split}
&
M(a)\in G\cL_2
\Longleftrightarrow
a\in G\big(K_{p,0}^{1/p+\gamma,0}\cap(c^\gamma+H^\infty)\big),
\\
&
M(a)\in G\cL_*
\Longleftrightarrow
a\in G\big(K_{p,q}^{1/p+\gamma,1/q-\gamma}\cap(c^\gamma+H^\infty)\big).
\end{split}
\end{equation}
Combining (\ref{eq:IA-1})--(\ref{eq:IA-3}), we get the assertion (b).
\end{proof}
Now we prove that the intersection of a generalized Krein
algebra $K_{p,q}^{\alpha,\beta}$ (or $K_{0,q}^{0,\beta}$) with
$c^\gamma+\overline{H^\infty}$ is inverse closed in $c^\gamma+\overline{H^\infty}$
(and thus in $C+\overline{H^\infty}$) if $\beta>1/q$.
\begin{lemma}\label{le:IAprime}
Let $0<\lambda\le 1$, $1<p,q<\infty$, $1/p+1/q=\lambda$, and $0<\gamma<\lambda-1/q$.
Suppose $K$ is either $K_{p,q}^{1/p-\gamma,1/q+\gamma}$ or $K_{0,q}^{0,1/q+\gamma}$.
\begin{enumerate}
\item[(a)]
The set $K\cap (c^\gamma+\overline{H^\infty})$ is a
Banach algebra under the quasi-sub\-multiplicative norm
\[
\|a\|:=
\|M(a)\|_{\cB(\ell_2^{-\gamma/2,\gamma/2})}+\|a\|_K.
\]

\item[(b)]
If $a\in K\cap (c^\gamma+\overline{H^\infty})$, then
\[
a\in G\big(K\cap(c^\gamma+\overline{H^\infty})\big)
\Longleftrightarrow
a\in G(c^\gamma+\overline{H^\infty}).
\]
\end{enumerate}
\end{lemma}
\begin{proof}
From Theorem~\ref{th:CH-algebra} and Theorem~\ref{th:M}(a) it follows that
$c^\gamma+\overline{H^\infty}$ is a Banach algebra under the Banach algebra norm
$\|a\|:=\|M(a)\|_{\cB(\ell_2^{-\gamma/2,\gamma/2})}$.
Hence we can finish the proof by making obvious modifications in the proof
of Lemma~\ref{le:IA}.
\end{proof}
We are in a position to prove our second main result. We will show that the
intersection of a generalized Krein algebra $K$ with $c^\gamma+H^\infty$ or
with $c^\gamma+\overline{H^\infty}$ is dense in $K$. This allows us to prove
that $K$ is inverse closed in $C+H^\infty$ or in $C+\overline{H^\infty}$,
respectively.
\begin{proof}[Proof of Theorem~{\rm\ref{th:invertibility}}]
(a) For $\alpha=1/p$ the statement of Theorem~\ref{th:invertibility}(a)
is nothing else than Theorem~\ref{th:Krein-invertibility}(b), (d).

Assume that $\alpha>1/p$ and put $\gamma:=\alpha-1/p$. Let $K$ be either
$K_{p,0}^{1/p+\gamma,0}$ or $K_{p,q}^{1/p+\gamma,1/q-\gamma}$ and
denote $B:=K\cap(c^\gamma+H^\infty)$. Then, by Lemma~\ref{le:embedding-Krein}(a),
$K\subset C+H^\infty$. In view of Lemma~\ref{le:IA}(a), $B$ is continuously
embedded in $K$.

Due to Lemma~\ref{le:embedding-Besov}, there exists a constant $c\in(0,\infty)$ such
that $\|f\|_{L^\infty}\le c\|f\|_{B_p^{1/p+\gamma}}$ for all $f\in B_p^{1/p+\gamma}$.
Let $a\in K$. By Lemma~\ref{le:density-Besov},  for every $\varepsilon>0$
there exists a Laurent polynomial $p$ such that
\[
\|Qa-p\|_{B_p^{1/p+\gamma}}<\frac{\varepsilon}{(c+1)\|Q\|_{\cB(B_p^{1/p+\gamma})}}.
\]
Hence $Qp+Pa\in B$ and
\[
\begin{split}
\|a-(Qp+Pa)\|_K &= \|Q(Qa-p)\|_K =
\|Q(Qa-p)\|_{L^\infty}+\|Q(Qa-p)\|_{B_p^{1/p+\gamma}}
\\
&\le(c+1)\|Q\|_{\cB(B_p^{1/p+\gamma})}\|Qa-p\|_{B_p^{1/p+\gamma}}<\varepsilon.
\end{split}
\]
Thus, $B$ is dense in $K$.

Combining Theorem~\ref{th:invertibility-cH} and Lemma~\ref{le:IA}(b), we obtain
that if $b\in B$, then the spectrum of $b$ in $B$ coincides with the spectrum
of $b$ in $C+H^\infty$. Hence for each $\varphi\in M_B$ there is a
$\psi\in M_{C+H^\infty}$ depending on $\varphi$ and $b$ such that
\begin{equation}\label{eq:invertibility-1}
\psi(b)=\varphi(b).
\end{equation}

Assume now that there exists an $a\in K$ being invertible in $C+H^\infty$
but not invertible in $K$. Thus, by Gelfand's theorem, we can find a
$\varphi\in M_K$ with $\varphi(a)=0$. Since $B$ is dense in $K$, there
exists a sequence $\{b_n\}\subset B$ such that
\begin{equation}\label{eq:invertibility-2}
\|a-b_n\|_K\to 0
\quad\mbox{as}\quad
n\to\infty.
\end{equation}
Consequently, $\varphi(b_n)\to\varphi(a)=0$ as $n\to\infty$. Since
$\varphi\in M_K\subset M_B$, it follows from (\ref{eq:invertibility-1})
that there are $\psi_n\in M_{C+H^\infty}$ such that
\begin{equation}\label{eq:invertibility-3}
\psi_n(b_n)=\varphi(b_n)\to 0
\quad\mbox{as}\quad
n\to\infty.
\end{equation}
On the other hand, from (\ref{eq:invertibility-2}) and the definition of the
norm in $K$ we get $\|a-b_n\|_{L^\infty}\to 0$ as $n\to\infty$. Since
$a\in G(C+H^\infty)$, it results that $b_n\in G(C+H^\infty)$ for all
sufficiently large $n$. Hence $\|b_n^{-1}-a^{-1}\|_{L^\infty}\to 0$
as $n\to\infty$. Thus $\|b_n^{-1}\|_{L^\infty}\le 2\|a^{-1}\|_{L^\infty}$
for all $n$ large enough. Since every multiplicative linear functional has norm
$1$, we obtain
\[
|\psi_n(b_n^{-1})|
\le
\|\psi_n\|\,\|b_n^{-1}\|_{L^\infty}
\le
2\|a^{-1}\|_{L^\infty}
\]
for all sufficiently large $n$. From this inequality and $\psi_n(b_n)\psi_n(b_n^{-1})=1$
it follows that $|\psi_n(b_n)|\ge (2\|a^{-1}\|_{L^\infty})^{-1}$ for all
$n$ large enough, but this contradicts (\ref{eq:invertibility-3}).
Part (a) is proved.

(b) The proof is analogous to the proof of part (a). We only need to replace
Lemma~\ref{le:IA} by Lemma~\ref{le:IAprime}.
\end{proof}
\section{Proof of higher order asymptotic formulas}
\label{sec:proof-asymptotic}
\subsection{Higher order asymptotic formulas: an abstract version}
Our asymptotic analysis is based on the following fact \cite[Section~10.34]{BS06}.
\begin{lemma}\label{le:asymptotics}
Suppose $a\in L_{N\times N}^\infty$ satisfies the following hypotheses:
\begin{enumerate}
\item[(i)]
there are two factorizations $a=u_-u_+=v_+v_-$, where
$u_-,v_-\in G(\overline{H^\infty})_{N\times N}$ and
$u_+,v_+\in G(H^\infty)_{N\times N}$;

\item[(ii)]
$u_-\in G(C+H^\infty)_{N\times N}$ or
$u_+\in G(C+\overline{H^\infty})_{N\times N}$.
\end{enumerate}
Then $D_n(a)\ne 0$ and
\begin{equation}\label{eq:asymptotics}
\frac{G(a)^{n+1}}{D_n(a)}=\det\Bigg(I-\sum_{k=0}^\infty F_{n,k}\Bigg)
\end{equation}
for sufficiently large $n$, where $G(a)$ is defined by {\rm (\ref{eq:def-G})}
and $F_{n,k}$ are defined by {\rm (\ref{eq:def-F})}.
\end{lemma}
From (i) and (ii) it follows that $a\in G(C+H^\infty)_{N\times N}$
or $a\in G(C+\overline{H^\infty})_{N\times N}$. Moreover, the operators $T(a)$
and $T(\widetilde{a})$ are invertible on $H_N^2$. Hence, in view of
Lemma~\ref{le:constant-G}, the constant $G(a)$ is well defined.

An earlier version of Lemma~\ref{le:asymptotics} was obtained in \cite{BS80}
(see also \cite[Section~6.15]{BS83}) with
\begin{enumerate}
\item[(iii)]
$u_-\in C_{N\times N}$ or $u_+\in C_{N\times N}$
\end{enumerate}
in place of (ii).

If the series $\sum_{k=0}^\infty F_{n,k}$ converges in the
topology of the Schatten-von Neumann class $\cC_m(H_N^2)$,
$m\in\N$, then one can remove its remainder in
(\ref{eq:asymptotics}) because the mapping
\[
\mathrm{det}_m(I+\cdot):\cC_m(H_N^2)\to\C
\]
is continuous (see, e.g. \cite[Theorem~6.5]{S77}).
More precisely, starting from (\ref{eq:asymptotics}) one can get the following.
\begin{theorem}\label{th:abstract}
Let $a\in L_{N\times N}^\infty$ satisfy the hypotheses {\rm (i)} and {\rm (ii)}
of Lemma~{\rm\ref{le:asymptotics}}. Define the constant $G(a)$, the functions
$b,c$, and the operators $F_{n,k}$ by {\rm (\ref{eq:def-G}), (\ref{eq:def-bc}),}
and {\rm (\ref{eq:def-F})}, respectively. If $m\in\N$ and $H(\widetilde{c})H(b)$,
$H(b)H(\widetilde{c})$ belong to the Schatten-von Neumann class $\cC_m(H_N^2)$,
then {\rm (\ref{eq:formula})} holds.
\end{theorem}
This statement is proved in \cite[Theorem~15]{K06} under the hypotheses (i)
and (iii). Exactly the same proof works under the hypotheses (i) and (ii)
because of Lemma~\ref{le:asymptotics}. Theorem~\ref{th:abstract}
is contained implicitly in \cite[Section~5]{BS80} (see also \cite[Theorem~6.20]{BS83},
\cite[Theorem~10.37]{BS06}).

Notice that we need the statement in this (new) form because
generalized Krein algebras contain discontinuous functions and therefore
the (old) hypothesis (iii) is not satisfied for generalized Krein algebras.
\subsection{Wiener-Hopf factorization in decomposing algebras}
Mark Krein \cite{Kr58} was the first to understand the Banach algebraic background
of Wiener-Hopf factorization and to present the method in a crystal-clear manner.
Gohberg and Krein \cite{GK58} proved that $a\in GW_{N\times N}$ admits a
Wiener-Hopf factorization. Later Budyanu and Gohberg developed an abstract
factorization theory in decomposing algebras of \textit{continuous} functions.
Their results are contained in \cite[Chap.~2]{CG81}. However, generalized Krein algebras
contain discontinuous functions and the results of Budyanu and Gohberg are not
applicable to them. One of the authors and Heinig \cite{HS84} extended
the theory of Budyanu and Gohberg to the case of decomposing algebras which may
contain \textit{discontinuous} functions. The following definitions and results are
taken from \cite{HS84} (see also \cite[Chap.~5]{BS83}).

Let $A$ be a Banach algebra of complex-valued functions on the unit circle
$\T$  under a Banach algebra norm $\|\cdot\|_A$. The algebra $A$ is said to be
\textit{decomposing} if it possesses the following properties:

\begin{enumerate}
\item[(a)]
$A$ is continuously embedded in $L^\infty$;

\item[(b)]
$A$ contains all Laurent polynomials;

\item[(c)]
$PA\subset A$ and $QA\subset A$.
\end{enumerate}

Using the closed graph theorem it is easy to deduce from (a)--(c)
that $P$ and $Q$ are bounded on $A$ and that $PA$ and $QA$ are
closed subalgebras of $A$. Given a decomposing algebra $A$ put
\[
A_+=PA,\quad
\stackrel{\circ}{A}_-=QA,\quad
\stackrel{\circ}{A}_+=\chi_1A_+,\quad
A_-=\chi_1\stackrel{\circ}{A}_-.
\]

Let $A$ be a decomposing algebra. A matrix function $a\in A_{N\times N}$
is said to \textit{admit a right Wiener-Hopf factorization in}
$A_{N\times N}$ if it can be represented in the form $a=a_-da_+$, where
$a_\pm\in G(A_\pm)_{N\times N}$ and
\[
d=\diag(\chi_{\kappa_1},\dots,\chi_{\kappa_N}),
\quad
\kappa_i\in\Z,
\quad
\kappa_1\le\kappa_2\le\dots\le\kappa_N.
\]
The integers $\kappa_i$ are usually called the \textit{right partial
indices} of $a$; they can be shown to be uniquely determined by $a$.
If $\kappa_1=\dots=\kappa_N=0$, then the Wiener-Hopf factorization is
said to be \textit{canonical}.
A decomposing algebra $A$ is said to have the \textit{factorization
property} if every matrix function in $GA_{N\times N}$ admits a right
Wiener-Hopf factorization in $A_{N\times N}$.

Let $\cR$ be the restriction to the unit circle $\T$ of the set of all
rational functions defined on the whole plane $\C$ and having no poles
on $\T$.
\begin{theorem}\label{th:HS}
Let $A$ be a decomposing algebra. If at least one of the sets
$(\cR\,\cap\stackrel{\circ}{A_-})+A_+$ or $\stackrel{\circ}{A_-}+(\cR\cap A_+)$
is dense in $A$, then $A$ has the factorization property.
\end{theorem}
\subsection{Fredholmness and invertibility of Toeplitz operators}
In this subsection we collect some well known facts about the Fredholmness
and invertibility of Toeplitz operators on $H_N^2$.
\begin{theorem}\label{th:Toeplitz-Fredholmness}
Let $a\in L_{N\times N}^\infty$.
\begin{enumerate}
\item[(a)]
If $a\in W_{N\times N}$, then $T(a)$ is Fredholm on $H_N^2$ if and only if
$\det a\in GW$.

\item[(b)]
If $a\in(C+\overline{H^\infty})_{N\times N}$, then $T(a)$ is Fredholm on $H_N^2$
if and only if $\det a$ belongs to $G(C+\overline{H^\infty})$.
\end{enumerate}
\end{theorem}
Part (a) is essentially due to Gohberg and Krein \cite{GK58}, part (b) is due
to Douglas, its proof is given in \cite[Theorem~2.94(a)]{BS06}.

The following result was proved by Widom and in a slightly different setting by
Simonenko.
\begin{theorem}\label{th:Toeplitz-invertibility}
If $a\in GL_{N\times N}^\infty$, then $T(a)$ is invertible on
$H_N^2$ if and only if the following two conditions hold:
\begin{enumerate}
\item[(a)]
$a$ admits a canonical right generalized factorization in $L_N^2$,
that is, there exist $a_-$, $a_+$ such that
\[
a=a_-a_+, \quad a_-^{\pm 1}\in(\overline{H^2})_{N\times N}, \quad
a_+^{\pm 1}\in(H^2)_{N\times N};
\]
\item[(b)]
the operator $M(a_-)PM(a_-^{-1})$ is bounded on $L_N^2$.
\end{enumerate}
\end{theorem}
For a proof, see, e.g. \cite[Chap.~7, Theorem~3.2]{CG81} or
\cite[Theorem~3.14]{LS87}.
\begin{theorem}\label{th:Toeplitz-tilde-invertibility}
If $a\in GL_{N\times N}^\infty$, then $T(\widetilde{a})$ is invertible on
$H_N^2$ if and only if $T(a^{-1})$ is invertible on $H_N^2$.
\end{theorem}
For a proof, see \cite[Proposition~7.19(b)]{BS06}.
\subsection{Wiener-Hopf factorization in generalized Krein algebras}
Let us show that all algebras mentioned in Theorem~\ref{th:main} have the
factorization property.
\begin{lemma}\label{le:fact-property}
Let $1<p,q<\infty$ and $0<\alpha,\beta<1$. Each of the algebras
\[
W\cap K_{p,0}^{1/p,0},
\quad
K_{p,0}^{\alpha,0} \mbox{ with } \alpha>1/p,
\quad
W\cap K_{0,q}^{0,1/q},
\quad
K_{0,q}^{0,\beta} \mbox{ with } \beta>1/q,
\quad
W\cap K_{p,q}^{1/p,1/q},
\]
and $K_{p,q}^{\alpha,\beta}$ with $\alpha\ne 1/p$, $1/p+1/q=\alpha+\beta\in(0,1]$
is a decomposing algebra with the factorization property.
\end{lemma}
\begin{proof}
This statement is proved for $W\cap K_{p,q}^{1/p,1/q}$ in
\cite[Section~10.24]{BS06}, the same argument applies also to
$W\cap K_{p,0}^{1/p,0}$ and $W\cap K_{0,q}^{0,1/q}$.

Under the imposed conditions on the parameters $p,q,\alpha,\beta$,
by Lemma~\ref{le:PQ}, the projections $P$ and $Q$ are bounded on each
$K_{p,0}^{\alpha,0}$, $K_{0,q}^{0,\beta}$, and $K_{p,q}^{\alpha,\beta}$.
Hence each of these algebras is decomposing.

If $\beta>1/q$, then, by Lemma~\ref{le:embedding-Besov}, $B_q^\beta$ is
continuously embedded in $C$ and, by Lemma~\ref{le:density-Besov}, the Laurent
polynomials are dense in $B_q^\beta$. It follows that $\cR\cap PB_q^\beta$
is dense in $PK_{p,q}^{\alpha,\beta}$ and in $PK_{0,q}^{0,\beta}$.
Theorem~\ref{th:HS} gives the factorization property of the algebras
$K_{p,q}^{\alpha,\beta}$ and $K_{0,q}^{0,\beta}$ if $\beta>1/q$. Analogously
one can show that $K_{p,q}^{\alpha,\beta}$ and $K_{p,0}^{\alpha,0}$ have the
factorization property if $\alpha>1/p$.
\end{proof}
We are ready to start the proof of our last main result.
\begin{proof}[Proof of Theorem~{\rm\ref{th:main}(a)}]
From Lemma~\ref{le:embedding-Krein} we conclude that
$K_{N\times N}\subset (C+H^\infty)_{N\times N}$ or
$K_{N\times N}\subset (C+\overline{H^\infty})_{N\times N}$.
In the first case $\widetilde{a}\in (C+\overline{H^\infty})_{N\times N}$.
Since $T(\widetilde{a})$ is invertible, from Theorem~\ref{th:Toeplitz-invertibility}(b)
we get $\det\widetilde{a}\in G(C+\overline{H^\infty})$ (and $\det\widetilde{a}\in GW$
if, in addition, $a\in W_{N\times N}$ due to Theorem~\ref{th:Toeplitz-invertibility}(a)).
Hence $\det a\in G(C+H^\infty)$. By Theorem~\ref{th:invertibility}, $\det a\in GK$.
Then, in view of \cite[Chap.~1, Theorem~1.1]{K87}, $a\in GK_{N\times N}$.
If $K_{N\times N}\subset(C+H^\infty)_{N\times N}$, then, as before, the
Fredholmness of $T(a)$ implies $a\in GK_{N\times N}$.

By Theorem~\ref{th:Toeplitz-invertibility}, there exist $a_-$, $a_+$ such that
$a=a_-a_+$, $a_-^{\pm 1}\in (\overline{H^2})_{N\times N}$,
and $a_-^{\pm 1}\in (H^2)_{N\times N}$. On the other hand, in view
of Lemma~\ref{le:fact-property}, $a$ admits a right Wiener-Hopf factorization
in $K_{N\times N}$, that is, there exist $u_\pm\in G(K_\pm)_{N\times N}$
such that $a=u_-du_+$. It is clear that
\[
u_-^{\pm 1}\in (K_-)_{N\times N}\subset(\overline{H^2})_{N\times N},
\quad
u_+^{\pm 1}\in (K_+)_{N\times N}\subset(H^2)_{N\times N}.
\]
Hence $a=u_-du_+$ is a right generalized factorization of $a$ in $L_N^2$.
It is well known that the set of partial indices of such a factorization
is unique (see, e.g. \cite[Corollary~2.1]{LS87}). Thus $d=1$ and $a=u_-u_+$.

In view of Theorem~\ref{th:Toeplitz-tilde-invertibility}, $T(a^{-1})$ is invertible
on $H_N^2$. By what has just been proved, there exist $f_\pm\in G(K_\pm)_{N\times N}$
such that $a^{-1}=f_-f_+$. Put $v_\pm:=f_\pm^{-1}$. In that case
$v_\pm\in G(K_\pm)_{N\times N}$
and $a=v_+v_-$.
\end{proof}
\subsection{Products of Hankel operators in Schatten-von Neumann classes}
Now we prove simple sufficient conditions guaranteeing the membership of
products of two Hankel operators in Schatten-von Neumann classes.
\begin{lemma}\label{le:product1}
Let $1<p,q<\infty$ and $0<\alpha,\beta<1$.
\begin{enumerate}
\item[(a)]
If $\alpha\ge 1/p$ and $b,c\in K_{p,0}^{\alpha,0}$, then the operators
$H(\widetilde{c})H(b)$ and $H(b)H(\widetilde{c})$ belong to the Schatten-von
Neumann class $\cC_p(H^2)$.
\item[(b)]
If $\beta\ge 1/q$ and $b,c\in K_{0,q}^{0,\beta}$, then the operators
$H(\widetilde{c})H(b)$ and $H(b)H(\widetilde{c})$ belong to the Schatten-von
Neumann class $\cC_q(H^2)$.
\end{enumerate}
\end{lemma}
\begin{proof}
(a) If $c\in K_{p,0}^{\alpha,0}$, then $Qa\in B_p^\alpha$. Since $\alpha\ge 1/p$,
we have $B_p^\alpha\subset B_p^{1/p}$. From (\ref{eq:Peller-0}) we conclude that
$H(\widetilde{c})\in\cC_p(H^2)$. Thus
$H(\widetilde{c})H(b), H(b)H(\widetilde{c})\in\cC_p(H^2)$.

(b) The proof is analogous.
\end{proof}
\begin{lemma}\label{le:product2}
Suppose $1<p,q,r<\infty$, $0<\alpha,\beta<1$, and
\begin{equation}\label{eq:product}
1/p+1/q=\alpha+\beta\in(0,1],
\quad
-1/2<\alpha-1/p<1/2.
\end{equation}
If $1/r=1/p+1/q$ and $b,c\in K_{p,q}^{\alpha,\beta}$, then the operators
$H(\widetilde{c})H(b)$ and $H(b)H(\widetilde{c})$ belong to the Schatten-von
Neumann class $\cC_r(H^2)$.
\end{lemma}
\begin{proof}
Put $\gamma:=\alpha-1/p=1/q-\beta$. If $b,c\in K_{p,q}^{\alpha,\beta}$, then
$Qb,Qc$ belong to $B_p^{1/p+\gamma}$ and $Pb,Pc$ are in $B_q^{1/q-\gamma}$. From
(\ref{eq:product}) we easily get
\[
\min\{\gamma,0\}>\max\{-1/2,-1/p\},
\quad
\max\{\gamma,0\}<\min\{1/2,1/q\}.
\]
Then, by Theorem~\ref{th:Peller2},
\[
QM(c)P\in\cC_p(\ell_2^{\gamma,0}),
\
QM(c)P\in\cC_p(\ell_2^{0,\gamma}),
\
PM(b)Q\in\cC_q(\ell_2^{\gamma,0}),
\
PM(b)Q\in\cC_q(\ell_2^{0,\gamma}).
\]
Since $1/r=1/p+1/q$, from Lemma~\ref{le:Schatten-properties}(c) we
deduce that
\[
PM(b)QM(c)P\in\cC_r(\ell_2^{\gamma,0}),
\quad
QM(c)PM(b)Q\in\cC_r(\ell_2^{0,\gamma}).
\]
Applying the flip operator to the second operator, we obtain
\[
JQM(c)PM(b)QJ\in\cC_r(\ell_2^{\gamma,0}).
\]
Considering the compressions of the above operators to $\ell_2(\Z_+)\sim H^2$,
we finally get $H(b)H(\widetilde{c}), H(\widetilde{c})H(b)\in\cC_r(H^2)$.
\end{proof}
Now we are ready to finish the proof of our last main result.
\begin{proof}[Proof of Theorem~{\rm\ref{th:main}(b)}]
From Theorem~\ref{th:main}(a) we know that $b,c\in K_{N\times N}$.
Let $b=(b_{ij})_{i,j=1}^N$ and $c=(c_{kl})_{k,l=1}^N$. Since
$b_{ij},c_{kl}\in K$, from Lemmas~\ref{le:product1} and
\ref{le:product2} we conclude that
$H(\widetilde{c_{kl}})H(b_{ij})$ and
$H(b_{ij})H(\widetilde{c_{kl}})$ belong to $\cC_{1/\lambda}(H^2)$
for all $i,j,k,l\in\{1,\dots,N\}$, where $\lambda$ is defined by
(\ref{eq:def-lambda}). Hence
$H(\widetilde{c})H(b),H(b)H(\widetilde{c})\in\cC_{1/\lambda}(H_N^2)$.
Since $m\ge 1/\lambda$, we finally get $H(\widetilde{c})H(b),
H(b)H(\widetilde{c})\in\cC_m(H_N^2)$ by
Lemma~\ref{le:Schatten-properties}(a).
\end{proof}
\begin{proof}[Proof of Theorem~{\rm\ref{th:main}(c)}]
It is only necessary to verify the hypotheses (i) and (ii) of
Lemma~\ref{le:asymptotics}. By Theorem~\ref{th:main}(a),
$u_-,v_-\in G(K\cap\overline{H^\infty})_{N\times N}$, $u_+,v_+\in
G(K\cap H^\infty)_{N\times N}$. Hence $u_-,v_-\in
G(\overline{H^\infty})_{N\times N}$, $u_+,v_+\in
G(H^\infty)_{N\times N}$, and $u_\pm\in GK_{N\times N}$. From this
fact and Theorem~\ref{th:invertibility} we see that $u_-\in
G(C+H^\infty)_{N\times N}$ or $u_+\in
G(C+\overline{H^\infty})_{N\times N}$. Thus, the hypotheses (i)
and (ii) are met. By Theorem~\ref{th:main}(b),
$H(\widetilde{c})H(b)$ and $H(b)H(\widetilde{c})$ belong to
$\cC_m(H_N^2)$. Applying Theorem~\ref{th:abstract}, we get
(\ref{eq:formula}).
\end{proof}

\end{document}